\documentclass[letterpaper]{amsart}

\newfont{\bbf}{msbm10} 
\newfont{\bbfs}{msbm7} 

\newcommand{\euO}{\mathfrak O}
\newcommand{\euP}{\mathfrak P}

\newcommand{\euM}{\mathfrak M}

\newcommand{\bX}{\relax\ifmmode\mathchoice{\mbox{\bbf X}}{\mbox{\bbf X}}
       {\mbox{\bbfs X}}{{X}}\else {\bbf X}\fi}
\newcommand{\bY}{\relax\ifmmode\mathchoice{\mbox{\bbf Y}}{\mbox{\bbf Y}}
       {\mbox{\bbfs Y}}{{Y}}\else {\bbf Y}\fi}

\newcommand{\bZ}{\relax\ifmmode\mathchoice{\mbox{\bbf Z}}{\mbox{\bbf Z}}
       {\mbox{\bbfs Z}}{{Z}}\else {\bbf Z}\fi}
\newcommand{\bN}{\relax\ifmmode\mathchoice{\mbox{\bbf N}}{\mbox{\bbf N}}
      {\mbox{\bbfs N}}{{N}}\else {\bbf N}\fi}
\newcommand{\bF}{\relax\ifmmode\mathchoice{\mbox{\bbf F}}{\mbox{\bbf F}}
       {\mbox{\bbfs F}}{{F}}\else {\bbf F}\fi}

\newtheorem{thm}{Theorem}[section]
\newtheorem{lem}[thm]{Lemma}

\theoremstyle{definition}

\theoremstyle{remark}
\newtheorem{rem}[thm]{Remark}

\numberwithin{equation}{section}

\usepackage{amsmath} 
\usepackage{amssymb}
\usepackage{eepic}
\usepackage{rotating}
\usepackage{eufrak}

\usepackage{latexsym}
\begin{document}
\title[Galois Structure]{The Galois Structure of Ambiguous
Ideals in Cyclic Extensions of Degree $8$} 
\author{G. Griffith Elder}
\address{Department of Mathematics, University of Nebraska at Omaha,
Omaha, Nebraska 68132-0243} 
\email{elder@unomaha.edu}
\date{October 6, 2002}
\begin{abstract}
In cyclic, degree $8$ extensions of algebraic number fields $N/K$,
ambiguous ideals in $N$ are canonical $\mathbb{Z}[C_8]$-modules.  Their
$\mathbb{Z}[C_8]$-structure is determined here.  It is described in terms of
indecomposable modules and determined by ramification invariants.
Although infinitely many indecomposable $\mathbb{Z}[C_8]$-modules are
available (classification by Yakovlev), only 23 appear.
\end{abstract}
\maketitle
\section{Introduction}

We are concerned with the interrelationship between two basic objects
in algebraic number theory: the {\em ring of integers} and the {\em
Galois group}.  In particular, we seek to understand the {\em effect}
of the Galois group upon the ring of integers.  At the same time, we
are also interested in the Galois action upon other fractional
ideals. So that the action may be similar, we restrict ourselves to
{\em ambiguous ideals} -- those that are mapped to themselves by the
Galois group.  The setting for our investigation is the family of
$C_8$-extensions.  This choice is guided by by a result of E. Noether
as well as results in Integral Representation Theory.

\noindent{\bf Noether's Normal Integral Basis Theorem}.
 A finite
Galois extension of number fields $N/K$ is said to be {\em at most
tamely ramified} (TAME) if the factorization of each prime ideal
$\euP_K$ (of $\euO_K$) in $\euO_N$ results in exponents (degrees of
ramification) that are relatively prime to the ideal $\euP_K$.  A {\em
normal integral basis} (NIB) is said to exist if there is an element
$\alpha\in\euO_N$ (in the ring of integers of $N$) whose conjugates,
$\{\sigma\alpha: \sigma\in \mbox{Gal}(N/K)\}$, provide a basis for
$\euO_N$ over $\euO_K$ (the integers in $K$).

Noether proved ${\rm NIB} \Rightarrow{\rm TAME}$; moreover, for local
number fields ${\rm NIB} \Leftrightarrow{\rm TAME}$, tying the Galois
module structure of the ring of integers to the arithmetic of the
extension \cite{emmy}.  This is a {\em nice} effect -- NIB means that
the integers are isomorphic to the group ring,
$\euO_K[\mbox{Gal}(N/K)]$. It is similar to the effect of the Galois
group on the field itself ({\em i.e.}  Normal Basis Theorem). The
impact of her result is two-fold: (1) We are encouraged to
localize. (2) We are directed away from tamely ramified extensions -- toward 
wildly ramified extensions and $p$-groups (See \cite{miyata}).

\noindent{\bf Integral Representation Theory} (Restricted to $p$-groups $G$).

\noindent{\em Classification of Modules}.
The number of 
indecomposable modules over a group ring $\mathbb{Z}[G]$ is, in general, infinite.
Only $\mathbb{Z}[C_p]$ and
$\mathbb{Z}[C_{p^2}]$ are of
{\em finite type}.  Still, among those of {\em infinite type}, there are two
whose classifications are somehow manageable.  These are the
ones of so--called {\em tame type} \cite{diet}:
$\mathbb{Z}[C_2\times C_2 ]$ (classification by L. A. Nazarova \cite{naza}) and
$\mathbb{Z}[C_8]$ (classification by A. V. Yakovlev \cite{jakov:2}).

\noindent{\em Unique Decomposition}.
The Krull--Schmidt--Azumaya Theorem does not, in general, hold:
although a module over a group ring will decompose into indecomposable
modules, this decomposition may not be unique. Fortunately, it does
hold for a few group rings, including $\mathbb{Z}[C_2\times C_2 ]$ and
$\mathbb{Z}[C_8]$ \cite{kling:ks}.

\noindent{\bf Topic}. Let $G=\mbox{Gal}(N/K)$.
We are led to ask the following
natural question: {\em What is the 
$\mathbb{Z}[G]$-module structure of ambiguous 
ideals when}
\begin{itemize}
\item {\em the number theory is `bad' (wild ramification), while}
\item {\em the representation theory is `good' (tame type, K--S--A)?}
\end{itemize}
In other words: {\em What is the $\mathbb{Z}[G]$-module structure of
ambiguous ideals in wildly ramified $C_2\times C_2$ and $C_8$ number
field extensions?}  Previous work solved this for $C_2\times
C_2$-extensions \cite{elder:biquad}, \cite{elder:byott}. So our focus
here is on $C_8$-extensions.  (Note: This question has already
been addressed for those group rings with `very good' representation
theory, those of {\em finite type}. See \cite{martha} and
\cite{elder:annals}.)

As with $C_2\times C_2$-extensions, the
$\mathbb{Z}[G]$-module structure of ambiguous ideals in
$C_8$-extensions is completely determined by the structure at its
$2$-adic completion -- our global question reduces to a
collection of local ones. We leave it to the reader to fill in the
details. (One may follow \cite[\S 2]{elder:biquad} using
\cite{wiegand}.)

\subsection{Local Question, Answer}
Let $K_0$ be a finite extension of the $2$-adic numbers $\mathbb{Q}_2$ and
let $K_n$ be a wildly ramified, cyclic, degree $2^n$ extension of
$K_0$ with $G=\mbox{Gal}(K_n/K_0)$.  The maximal ideal $\euP_n$
in $K_n$ is unique (therefore ambiguous). So every fractional ideal
$\euP_n^i$ is ambiguous.  We ask: {\em What is the
$\mathbb{Z}_2[G]$-module structure of $\euP_n^i$ for $n=1, 2, 3$?} ($\mathbb{Z}_2$
denotes the $2$-adic integers.)  The answered is given by the following
theorem and the description of the modules $\euM_s(i,b_1, \cdots , b_s)$.

Let $T$ denote the maximal unramified extension of $\mathbb{Q}_2$ in
$K_0$. Following \cite[Ch IV]{serre:local}, let $G=G_{-1}\supseteq G_0
\supseteq G_1 \supseteq \cdots$ denote the ramification filtration.
Use subscripts to denote field of reference, so $\euO_k$ denotes the
ring of integers of $k$.
\begin{thm}
Let $K_n/K_0$ be a cyclic extension of degree $2^n$ and let
$k\subseteq K_0$ be an unramified extension of $\mathbb{Q}_2$.  
Suppose that $|G_1|\leq 8$ {\rm (}i.e. $s=1, 2$ or $3${\rm )} and
let $b_1,
\cdots , b_s$ be the break numbers in the ramification filtration of
$G_1$. 
If $\euM_s(i,b_1, \cdots , b_s)$ is the $\mathbb{Z}_2[G_1]$-module defined
below, then 
$$\euP_n^i\cong\euO_k[G]\otimes_{\mathbb{Z}_2[G_1]} \euM_s(i,b_1,
\cdots , b_s)^{[T:k]}\mbox{ as left $\euO_k[G]$-modules.}$$ 
\end{thm}

\subsubsection{$\euM_s(i,b_1, \cdots , b_s)$}
Indecomposable
modules are listed in Appendix A and $e_0$ denotes the absolute ramification index of $K_0$. 
Following \cite{martha} and \cite{elder:annals},
\begin{equation}
\euM_1(i,b_1)=
(\mathcal{R}_0\oplus\mathcal{R}_1)^{\lceil (i+b_1)/2\rceil -\lceil i/2\rceil}
\oplus
\mathbb{Z}_2[G_1]^{e_0-(\lceil (i+b_1)/2\rceil -\lceil i/2\rceil)}
\end{equation}

\begin{equation}
\euM_2(i,b_1,b_2)=\mathcal{I}^{a}\oplus
\begin{cases} 
\mathcal{H}^{b_A}\oplus
\mathcal{G}^{c_A}\oplus
\mathcal{L}^{d_A}
& \text{if $b_2+2b_1>4e_0$}, \\
\mathcal{H}^{b_B}\oplus
\mathcal{M}^{c_B}\oplus
\mathcal{L}^{d_B} 
&  \text{if $b_2+2b_1<4e_0$}.
\end{cases}
\end{equation}
where 
$a=\lceil (i+b_2)/4\rceil - \lceil (i+2b_1)/4\rceil$,
$b_A=e_0+\lceil i/4\rceil - \lceil (i+b_2)/4\rceil$,
$b_B=\lceil (i+b_2+2b_1)/4\rceil - \lceil (i+b_2)/4\rceil$,
$c_A=\lceil (i+b_2+2b_1)/4\rceil - e_0- \lceil i/4\rceil$,
$c_B=e_0+\lceil i/4\rceil - \lceil (i+b_2+2b_1)/4\rceil$,
$d_A=e_0+\lceil (i+2b_1)/4\rceil - \lceil (i+b_2+2b_1)/4\rceil$,
$d_B=\lceil (i+2b_1)/4\rceil - \lceil i/4\rceil$.

The description of $\euM_3(i,b_1,b_2,b_3)$ is given by Tables 1 and 2.
Note the eight columns in each table.  There are eight cases. 
Each module that appears in
$\euM_3(i,b_1,b_2,b_3)$, except for $\mathcal{R}_3$, is listed in the
appropriate column of Table 1. The multiplicity of the module is
appears in the corresponding spot in Table 2. The multiplicity of
$\mathcal{R}_3$ follows the tables.

\renewcommand{\arraystretch}{.8}
\setlength{\tabcolsep}{2.5mm}
\begin{sideways}
\begin{tabular}{l}
{\bf Table 1.} \vspace*{1mm}\\
\begin{tabular}{cccccccc} \vspace*{1mm}
$ {\mathbf A} $&$ {\mathbf B} $&$ {\mathbf C} $&$ {\mathbf D} $&$ {\mathbf E} $&$ {\mathbf F} $&$ {\mathbf G} $&$ {\mathbf H}  $\\  
$\mathcal{H} $&$ \mathcal{H} $&$ \mathcal{H} $&$ \mathcal{H} $&$ \mathcal{I}_1 $&$\mathcal{I}_1 $&$\mathcal{I}_1 $&$\mathcal{I}_1 $\\
$ $&$  $&$ \mathcal{H}_1\mathcal{L} $&$ \mathcal{H}_1\mathcal{L}  $&$ \mathcal{H}_1\mathcal{L} $&$ \mathcal{H}_1\mathcal{L}  $&$ $&$ $\\
$ $&$  $&$  $&$ \mathcal{H}_1\mathcal{G}  $&$  $&$ \mathcal{H}_1\mathcal{G}  $&$ $&$ $\\
$\mathcal{H}_2 $&$ \mathcal{H}_2 $&$ \mathcal{H}_1  $&$ \mathcal{H}_1  $&$ \mathcal{H}_1 $&$ \mathcal{H}_1  $&$\mathcal{H}_1 $&$\mathcal{H}_1 $\\
$ \mathcal{M} $&$ \mathcal{H}_{1,2} $&$ \mathcal{H}_{1,2}  $&$\mathcal{H}_{1,2} $&$\mathcal{G} $&$\mathcal{G} $&$\mathcal{G} $&$\mathcal{G} $\\
$ \mathcal{M}_1 $&$ \mathcal{M}_1 $&$ \mathcal{G}_4 $&$ \mathcal{G}_4  $&$\mathcal{G}_4 $&$\mathcal{G}_4 $&$\mathcal{G}_4 $&$\mathcal{G}_4 $\\
$\mathcal{L} $&$ \mathcal{L} $&$ \mathcal{L} $&$ \mathcal{L}  $&$\mathcal{G}_3 $&$\mathcal{G}_3 $&$\mathcal{G}_3 $&$\mathcal{G}_3 $\\
$ \mathcal{L}_3 $&$ \mathcal{L}_3 $&$ \mathcal{L}_3 $&$ \mathcal{L}_3  $&$\mathcal{L}_3 $&$\mathcal{L}_3 $&$\mathcal{G}_2 $&$\mathcal{G}_2 $\\
$ \mathcal{I} $&$ \mathcal{I} $&$ \mathcal{I} $&$ \mathcal{L}_2  $&$\mathcal{I} $&$\mathcal{L}_2 $&$\mathcal{L}_2 $&$\mathcal{G}_1 $\\
$ \mathcal{I}_2 $&$ \mathcal{I}_2 $&$ \mathcal{I}_2 $&$ \mathcal{I}_2  $&$\mathcal{I}_2 $&$\mathcal{I}_2 $&$\mathcal{L}_1 $&$\mathcal{L}_1 $
\end{tabular}\vspace*{10mm}\\
{\bf Table 2.} \vspace*{1mm}\\
\begin{tabular}{cccccccc}  \vspace*{1mm}
$ {\mathbf A} $&$ {\mathbf B} $&$ {\mathbf C} $&$ {\mathbf D} $ & ${\mathbf E} $&$ {\mathbf F} $&$ {\mathbf G} $&$ {\mathbf H}  $\\  
$\bar{d}-b  $&$ \bar{d}-b $&$ \bar{a}-b-e_0 $&$ \bar{a}-b-e_0 $ &$ b+e_0-\bar{a} $&$ b+e_0-\bar{a} $&$b-c $&$b-c  $\\ 
$ $&$  $&$ w+e_0-\bar{a} $&$ \bar{y}+m-\bar{a}  $ &$ w-b $&$ \bar{y}+m-e_0-b  $&$ $&$ $\\
$ $&$  $&$  $&$ d-\bar{y}+e_0  $ &$  $&$ \bar{c}-\bar{y}  $&$ $&$ $\\
$d+e_0-\bar{d} $&$ \bar{a}-\bar{d}-e_0 $&$ \bar{d}-w $&$ \bar{d}-d-m  $& $ a-w $&$ a+e_0-\bar{c}-m  $&$a-b $&$a-b $\\
$ \bar{a}-d-2e_0 $&$ d+2e_0-\bar{a} $&$ a-\bar{d}  $&$a-\bar{d} $ &$\bar{d}-a $&$\bar{d}-a $&$\bar{d}-a $&$\bar{d}-a $\\
$ a+e_0-\bar{a} $&$ a-d-e_0 $&$ d+e_0-a $&$ \bar{w}-m-a  $&$
\bar{c}-\bar{d} $&$\bar{y}-\bar{d} $&$\bar{c}-\bar{d} $&$\bar{c}-\bar{d} $\\
$\bar{c}-a $&$ \bar{c}-a $&$ \bar{c}-d-e_0 $&$ \bar{c}-d-e_0  $&$ 
d+e_0-\bar{c} $&$d+e_0-\bar{c} $&$\bar{b}-\bar{c} $&$\bar{b}-\bar{c} $\\
$ c+e_0-\bar{c} $&$  c+e_0-\bar{c} $&$ \bar{z}+b_1-\bar{c}  $&$ \bar{z}+b_1-\bar{c}  $&$
y-d $&$y+e_0-d $&$ d+e_0-\bar{b} $&$\bar{a}-\bar{b} $\\
$ \bar{b}-c-e_0 $&$ \bar{b}-c-e_0 $&$ \bar{b}-c-e_0 $&$ c+e_0-\bar{b}  $&$ \bar{b}-c-e_0 $&$ c+e_0-\bar{b} $&$ d+e_0-\bar{b} $&$ d+e_0-\bar{a}  $\\
$ b-\bar{b}+e_0 $&$ b-\bar{b}+e_0 $&$ b-\bar{b}+e_0 $&$ b-c  $ & $\bar{a}-\bar{b} $&$\bar{a}-c-e_0 $&$c+e_0-\bar{a} $&$ c-d $ \\
\end{tabular}
\end{tabular}
\end{sideways}

\begin{itemize}
\item[{\em Cases}.]
\item[$A$.] $4e_0-4b_1/3<b_2$ (including {\em Stable Ramification}, $b_1\geq e_0$).
\item[$B$.] $4e_0-2b_1< b_2 < 4e_0-4b_1/3$
\item[$C$.] $4e_0-4b_1< b_2 < 4e_0-2b_1$ and $b_2>(4e_0+4b_1)/3$
\item[$D$.] $4e_0-4b_1< b_2 < 4e_0-2b_1$ and $b_2<(4e_0+4b_1)/3$
\item[$E$.] $b_2 < 4e_0-4b_1$ and $b_3>8e_0+4b_1-2b_2$
\item[$F$.] $b_2 < 4e_0-4b_1$ and $8e_0+4b_1-2b_2<b_3<8e_0+4b_1-2b_2$
\item[$G$.] $8e_0-4b_1-2b_2 < b_3 < 8e_0-2b_2$
\item[$H$.] $b_3<8e_0-4b_1-2b_2$
\end{itemize}
A graphic representation of these cases
appears in \S 3.2.
\vspace{1mm}

\noindent {\em Constants used in Table 2}.
$a:=\lceil (i-2b_2)/8\rceil$,
$\bar{a}:=\lceil (i+b_3-2b_2)/8\rceil$,
$b:=\lceil (i-2b_2-4b_1)/8\rceil$,
$\bar{b}:=\lceil (i+b_3-2b_2-4b_1)/8\rceil$,
$c:=\lceil(i-4b_2)/8\rceil$,
$\bar{c}:=\lceil(i+b_3-4b_2)/8\rceil$,
$d:=\lceil(i-4b_2-4b_1)/8\rceil$,
$\bar{d}:=\lceil(i+b_3-4b_2-4b_1)/8\rceil$,
$w:=\lceil(i-2b_2-2b_1)/8\rceil$,
$\bar{w}:=\lceil(i+b_3-2b_2-2b_1)/8\rceil$,
$y:=\lceil(i-4b_2-2b_1)/8\rceil$,
$\bar{y}:=\lceil(i+b_3-4b_2-2b_1)/8\rceil$,
$\bar{z}:=\lceil(i+b_3-4b_2-6b_1)/8\rceil$,
$m:=(b_2-b_1)/2$.
\vspace{1mm}

\noindent {\em The multiplicity of $\mathcal{R}_3$}.
In Cases $A$ and
$B$ it is
$((\bar{a}+\bar{b}+\bar{c}+\bar{d})-(a+b+c+d)-3e_0)f$.  In Case $C$, it is
$((\bar{a}+\bar{b}+\bar{c}+\bar{d})-(a+b+d)-2e_0-(\bar{z}+b_1))f$.
While in Case $D$ it is
$((\bar{a}+\bar{b}+\bar{c}+\bar{d})-(a+b)-e_0+m-(\bar{w}+\bar{z}+b_1))f$.
In Case $E$, it is $(\bar{b}+\bar{d}-a-y-e_0)f$.  In Case $F$ it is
$((\bar{b}+\bar{d}+\bar{c})-(a+y+\bar{y})-e_0)f$.  Finally, in Cases
$G$ and $H$, the number of $\mathcal{R}_3$ that appear
is $(\bar{d}-a)f$.

\subsection{Discussion}
$\mbox{ }$
\vspace{1mm}

\noindent{\bf Cyclic $p$-Extensions}.  The Galois module structure of
the ring of integers in fully and wildly ramified, cyclic, local
extensions of degree $p^n$ was studied in \cite{elder:jnt} and more
recently in \cite{elder:bord}. Both of these papers required a lower
bound on the first ramification number $b_1$. In particular,
\cite{elder:bord} restricted $b_1$ to about half of its possible
values, under so-called {\em strong ramification}.  In this paper, by
focusing on $p=2$ we are able to remove this restriction. Our
work sheds light (1) on {\em strong ramification} and (2) on the
 structures that are possible outside of it.

\noindent (1) {\em Strong ramification} for $p=2$ means
$b_1>e_0$, a small part of Case $A$.  The structure under
{\em strong ramification} given by \cite[Thm 5.3]{elder:bord}, when
restricted to $p=2$, remains valid throughout Case $A$. 
{\em What then should Case $A$ be, for odd $p$?}

\noindent (2) Suppose that `nice' refers to the structure under {\em
strong ramification}, indeed under Case $A$. Does the structure remain
relatively `nice' beyond Case $A$? This depends upon a
precise definition. Let an indecomposable module be {\em nice} if it
is made up of distinct irreducible modules. Note only {\em nice}
modules appear in Case $A$.  But then, as we leave Case $A$, the
structure turns {\em nasty} immediately.  At least one of
$\mathcal{H}_{1,2}$, $\mathcal{H}_1\mathcal{L}$ and
$\mathcal{H}_1\mathcal{G}$ appears in every Case $B$ through $F$.
\vspace{1mm}

\noindent{\bf Induced Structure}.  The subfield of $K_n$ fixed by the
first ramification group $G_1$ is tame over the base field $K_0$.
Miyata generalized Noether's Theorem proving that each ideal is {\em
relatively projective} over $G_1$ \cite{miyata}.  In other words, the
ideals are direct summands of modules that have been induced from
$G_1$ to $G$ \cite[\S10]{curt}.  We find, in our situation, that
ideals are {\em relatively free} over $G_1$. See \cite[Thm
2]{miyata:2} for a more general, related result.
\vspace{1mm}

\noindent{\bf Extension of Ground Ring}.
When studying the structure of ideals in an extension $K_n/K_0$ over a
group ring, one must choose a ring of coefficients. Does one study
`fine' structure -- over $\euO_0[G]$ where the coefficients are
integers in $K_0$. Does one study `coarse' structure -- over $\mathbb{Z}_2[G]$.
We study a canonical intermediate structure -- over $\euO_T[G]$ where
the coefficients belong to the Witt ring of the residue class field.
We determine this structure by listing generators and
relations. Interestingly, the coefficients in these
relations always belong $\mathbb{Z}_2[G]$. Therefore
$\euO_T[G]$-structure results, by extension of the ground ring, from
$\mathbb{Z}_2[G]$-structure \cite[\S30B]{curt}.
\vspace{1mm}

\noindent{\bf Realizable Modules}.
Let $\mathcal{S}_G$ denote the set of realizable indecomposable
$\mathbb{Z}_2[G]$-modules: Those indecomposable $\mathbb{Z}_2[G]$-modules that
appear in the decomposition of some ambiguous ideal in an extension
$N/K$ with $\mbox{Gal}(N/K)\cong G$.  Chinburg asked whether
$\mathcal{S}_G$ could be infinite.  In \cite{elder:biquad}, since
$\mathcal{S}_{C_2\times C_2}$ is
infinite, the answer
was found to be {\em yes}. We determine here that although the
set of indecomposable $\mathbb{Z}_2[C_8]$-modules is infinite, 
$\mathcal{S}_{C_8}$ is finite.
The sequence $|\mathcal{S}_{C_{2^n}}|$, $n=0, 1, 2, \ldots$
begins
$$ 1, 3, 7, 23\ldots$$

\subsection{Organization of Paper}
Preliminary results are presented in \S 2, main results in \S
3. There are two appendixes.  Appendix A lists all necessary
indecomposable modules.  Appendix B lists bases for our ideals.

\noindent Preliminary Results: In \S 2.1 we handle the special case
when a ramification break number is even. In \S 2.2, we present a
strategy for handling odd ramification numbers.  To motivate our work
in \S 3, we implement this strategy for $|G_1|=2$ and $4$, in \S 2.2.1
and \S 2.2.3 respectively.  We conclude, in \S 2.3, with a reduction
to totally ramified extensions.

\noindent Main Results: We begin in \S 3.1 with a brief outline and
discussion. Then, we catalog ramification numbers and prove some
technical lemmas in \S 3.2. All this sets the stage for our work in \S
3.3 determining the Galois structure of ideals in fully, though {\em
unstably}, ramified $C_8$-extensions. This is our primary focus. Our
work in \S 3.4 on {\em stably ramified} extensions is essentially
contained in \cite{elder:bord}.

\section{Preliminary Results}
We continue to use the notation of \S 1.1. Let $K_0$ be a finite
extension of $\mathbb{Q}_2$ and $K_n/K_0$ be a cyclic extension of degree
$2^n$. Let $\sigma$ generate $G=\mbox{Gal}(K_n/K_0)$ and use
subscripts to distinguish among subfields. So $K_i$ denotes the fixed
field of $\langle\sigma ^{2^i}\rangle$, $\euO_i$ denotes the ring of
integers of $K_i$ and $\euP_i$ denotes the maximal ideal of
$\euO_i$. Let $v_i$ be the additive valuation in $K_i$, $\pi_i$ its
prime element, so that $v_i(\pi_i)=1$.  Let $\mbox{Tr}_{i,j}$ denote
the trace from $K_i$ down to $K_j$.  Recall the ramification
filtration of $G$. Note $G_{-1}= G_0$ if and only if $K_n/K_0$ is
fully ramified. Also since $G$ is a $2$-group and $[G_1:G_0]$ is odd,
$G_0=G_1$. Furthermore since $G$ is cyclic and $G_i/G_{i-1}$ is
elementary abelian for $i>1$, there are $s=\log_2|G_1|$ breaks in the
filtration of $G_1$ \cite[p 67]{serre:local}.  Let $b_1<b_2<
\ldots<b_s$ denote these break numbers.  (The break numbers of $G$ may
include $-1$ as well.)  It is a standard exercise to show that $b_1,
\ldots,b_s$ are {\em all} either odd or even \cite[Ex 3, p 71]{serre:local}.
When they are even, we are in an extreme case, called {\em maximal
ramification}.  The general case, when they are odd, will be our
primary concern.
\subsection{Even Ramification Numbers}
If $b_1, \ldots,b_s$ are even, we use idempotent elements of the group
algebra, $\mathbb{Q}_2[G]$, and Ullom's generalization of Noether's result
\cite[Thm 1]{ullom}, to determine the structure of each ideal.  In
doing so, we rely upon two observations: (1) Idempotent elements in
$\mathbb{Q}_2[G]$ that map an ideal into itself, decompose the ideal. (2)
Modules over a principal ideal domain are free.

We illustrate this process in one case, leaving other cases to the
reader.  Suppose $|G|=8$ and $|G_1|=4$. So $K_3/K_0$ is only
partially ramified and $s=2$.  From \cite[IV \S 2 Ex 3]{serre:local},
$b_1=2e_0$ and $b_2=4e_0$.  Using \cite[V \S 3]{serre:local}, one
finds that $(1/2)(\sigma^4+1)\euP_3^i\subseteq \euP_3^i$. As a result,
the idempotent $(\sigma^4+1)/2$ decomposes the ideal $\euP_3^i\cong
\euP_2^{\lceil i/2\rceil}\oplus M_2$ with $(\sigma^4+1)M_2=0$.
Meanwhile $(1/2)(\sigma^2+1)\euP_2^{\lceil i/2\rceil} \subseteq
\euP_2^{\lceil i/2\rceil}$. So $\euP_2^{\lceil i/2\rceil}$ decomposes
as well. This yields $\euP_3^i\cong \euP_1^{\lceil i/4\rceil}\oplus
M_1\oplus M_2$ with $(\sigma^4+1)M_2=0$ and $(\sigma^2+1)M_1=0$.
Each $M_i$ may be viewed as a module over
$\euO_{T_K}[\sigma ]/(\sigma^{2^i}+1)$, a principal ideal domain. 
So
$M_i$ is free over $\euO_{T_K}[\sigma ]/(\sigma^{2^i}+1)$.  
Ullom's result provides a
normal integral basis for $\euP_1^{\lceil i/4\rceil}$.
Counting $\euO_T$-ranks, we find that
$$\euP_3^i\cong \frac{\euO_{T}[\sigma]}{(\sigma^2-1)}^{e_0}\oplus
\frac{\euO_{T}[\sigma]}{(\sigma^2+1)}^{e_0}\oplus
\frac{\euO_{T}[\sigma]}{(\sigma^4+1)}^{e_0}.$$

\subsection{Odd Ramification Numbers}Henceforth the ramification numbers 
will be odd. In this context we will
use the following
technical result (with
$K_i/K_{i-1}$).
\begin{lem} 
Let $k$ be a finite extension of $\mathbb{Q}_2$ and $K/k$ be a ramified
quadratic extension.  Let $e_k$ be the absolute ramification index of
$k$.  Assume that $\sigma$ generates the Galois group and that the
ramification number, $b<2e_k$, is odd.  Then
\begin{enumerate}
\item $v_K((\sigma\pm 1)\alpha)=v_K(\alpha)+b$ for $v_K(\alpha)$ odd;
\item if $\tau\in k$, there is a $\rho\in K$ such that $(\sigma
+1)\rho =\tau$ and $v_K(\rho )=v_K(\tau)-b$;
\item if $v_K(\alpha)$ is even and $(\sigma+1)\alpha=0$, there is a
$\theta\in K$ such that $\alpha=(\sigma -1)\theta$ and
$v_K(\theta)=v_K(\alpha)-b$.
\end{enumerate}
\end{lem} 
\begin{proof}
These may be shown using \cite[V \S3]{serre:local}, as in \cite[Lem
3.12--14]{elder:biquad}.
\end{proof}

\noindent Our {\em strategy} is based upon the following observations:
\newline (1) Under wild
ramification, Galois action `shifts/increases' valuation
(Lemma 2.1(1)). So an element may be used to `construct'
other elements with distinct valuation.  
\newline (2) Elements with distinct valuation may be used to construct
bases.  If the valuation map $v_n: K_n\rightarrow \mathbb{Z}$ is one--to--one
on a subset $A\subseteq K_n$, while $v_n(A)$ is onto $\{i, i+1, \ldots
, i+v_n(2)-1\}$; then $A$ is a basis for $\euP_n^i$ over the integers
in the maximal unramified subfield of $K_n$. If $K_n/K_0$ is fully
ramified, this subfield is $T$.
\newline The {\em strategy} is illustrated below. It is:
Use {\em Galois Action} to {\em Create Bases}.
\subsubsection{First Ramification Group of Order Two}
Suppose that $|G_1|=2$.  To use {\em Observation} (1), we pick
$\alpha\in K_n$ an element with $v_n(\alpha)=b_1$ ({\em e.g.}
$\alpha=\pi_n^{b_1}$).  Let $\alpha_m:=\alpha\cdot\pi_0^m$.  Since
$v_n(\pi_0)=2$, $v_n(\alpha_m)=b_1+2m$.  Use Lemma 2.1 with
$K_n/K_{n-1}$. So $v_n((\sigma^{2^{n-1}} +1)\alpha_m)=2b_1+2m$.  Since
$b_1$ is odd, the valuations of $\alpha_m$ and $(\sigma^{2^{n-1}}
+1)\alpha_m$ have opposite parity.  The valuations for all $m$ lie in
one--to--one correspondence with $\mathbb{Z}$. Select those with valuation in
$\{i, \ldots, i+v_n(2)-1\}$. Replace $\pi_0^{e_0}$ by $2$ whenever
possible. The result is
\begin{multline}
\mathcal{B}:=\left\{\alpha_m, (\sigma^{2^{n-1}} +1)\alpha_m :
\lceil (i-b_1)/2\rceil \leq m\leq e_0+
\lceil i/2\rceil -b_1-1 \right \} \\ \cup
\left\{(\sigma^{2^{n-1}} +1)\alpha_m, 2\alpha_m :
\lceil i/2\rceil -b_1 \leq m\leq
\lceil (i-b_1)/2\rceil -1 \right \}.
\end{multline}

Since $K_{n-1}/K_0$ is unramified, there is a root of unity $\zeta$
with $K_{n-1}=K_0(\zeta)$.  The maximal unramified extension $\mathbb{Q}_2$
in $K_n$ is $T(\zeta)$.  By {\em Observation} (2), $\mathcal{B}$
is a basis for $\euP_n^i$ over $\euO_{T(\zeta)}$.  Note that
$\euO_{T(\zeta)}\cdot\alpha_m + \euO_{T(\zeta)}\cdot(\sigma^{2^{n-1}}
+1)\alpha_m= \euO_{T(\zeta)}\cdot\alpha_m +
\euO_{T(\zeta)}\cdot\sigma\alpha_m$ yields the group ring,
$\euO_{T(\zeta)}[G_1]$, while $\euO_{T(\zeta)}\cdot(\sigma^{2^{n-1}}
+1)\alpha_m + \euO_{T(\zeta)}\cdot2\alpha_m =
\euO_{T(\zeta)}\cdot(\sigma^{2^{n-1}} +1)\alpha_m +
\euO_{T(\zeta)}\cdot(\sigma^{2^{n-1}} -1)\alpha_m$ yields the maximal
order of $\euO_{T(\zeta)}[G_1]$.  Restricting coefficients and
counting leads to the $\euO_k[G_1]$-module structure of $\euP_n^i$,
and to $\mathcal{M}_1(i,b_1)$ as in (1.1). 

Next, we extend $\mathcal{B}$ to a basis upon which the
action of the whole group can be followed.  Since $K_{n-1}/K_0$ is
unramified, there is a normal field basis for
$\euP_{n-1}^j/\euP_{n-1}^{j+1}$ over $\euO_{0}/\euP_{0}$ (for each
$j$).  Of course, $[\euO_0/\euP_0:\euO_T/\euP_T]=1$.  So
$\euP_{n-1}^j/\euP_{n-1}^{j+1}$ has a normal field basis over
$\euO_{T}/\euP_{T}$.  For $j=b_1$, this means that there is an element
$\mu\in\euP_{n-1}^{b_1}$ and basis $\mu, \sigma\mu, \ldots ,
\sigma^{2^{n-1}-1}\mu$.  Using Lemma 2.1(2), there is an 
$\alpha\in K_n$ with $v_n(\alpha)=b_1$ such that $(\sigma^{2^{n-1}}
+1)\alpha=\mu$.  Then $\alpha, \sigma\alpha, \ldots ,
\sigma^{2^{n-1}-1}\alpha$ is a normal field basis for
$\euP_n^{b_1}/\euP_n^{b_1+1}$ over $\euO_T/\euP_T$.  Since
$\{\sigma^j(\sigma^{2^{n-1}}+1 )\alpha:j=0, \ldots ,2^{n-1}-1\}$ is a
basis for $\euP_{n-1}^{b_1}/\euP_{n-1}^{b_1+1}$, it is also a basis
for $\euP_{n}^{2b_1}/\euP_{n}^{2b_1+1}$ over $\euO_T/\euP_T$. This
together with the fact that $\{\sigma^j\alpha: j=0, \ldots
,2^{n-1}-1\}$ is a basis for $\euP_{n}^{b_1}/\euP_{n}^{b_1+1}$ over
$\euO_T/\euP_T$ leads to 
$\cup_{j=0}^{2^{n-1}-1}\sigma^j\mathcal{B}$ being a basis for $\euP_n^i$
over $\euO_T$, and
$\euP_n^i\cong \euO_T[G]\otimes_{\mathbb{Z}_2[G_1]}\mathcal{M}_1(i,b_1)$
as
$\euO_T[G]$-modules.
\subsubsection{An Application of Nakayama's Lemma}
In the previous section we were able to follow the Galois action from
one basis element to another {\em explicitly}.  This level of detail
becomes overwhelming as we generalize to $|G_1|=4, 8$. Fortunately,
Nakayama's Lemma allows us to push some of these details into the
background.
\begin{lem}Let
$\mathcal{A}$ be a $\euO_k[C_{2^n}]$-module (torsion-free over
$\euO_k$) where $C_{2^n}=\langle \sigma\rangle$ and $k$ denotes an
unramified extension of $\mathbb{Q}_2$. Let $H$ denote the subgroup of order
$2$, $A^H$ the submodule fixed by $H$, and $\mbox{\rm Tr}_{H}A$ the
image under the trace. Then $\mbox{\rm Tr}_{H}\mathcal{A}/\left
((\sigma-1)\mbox{\rm Tr}_{H}\mathcal{A}+2\mathcal{A}^{H}\right )$ is
free over $\euO_k/2\euO_k$.  Suppose that
$\mathcal{B}\subseteq\mathcal{A}$ such that $\mbox{\rm
Tr}_{H}\mathcal{B}$ is a basis for $\mbox{\rm Tr}_{H}\mathcal{A}/\left
((\sigma-1)\mbox{\rm Tr}_{H}\mathcal{A}+2\mathcal{A}^{H}\right )$ then
$\mathcal{B}$ can be extended to a $\euO_k[C_{2^n}]/\langle \mbox{\rm
Tr}_H\rangle$-basis of $\mathcal{A}/\mathcal{A}^{H}$.
\end{lem}
\begin{proof} 
Since $\mathcal{A}/\mathcal{A}^{H}$ is a module over the principal
ideal domain $\euO_k[C_{2^n}]/\langle \mbox{Tr}_H\rangle$, it is free.
So $\mathcal{C}:=\mathcal{A}/\mathcal{A}^{H}\cong\left
(\euO_k[C_{2^n}]/\langle \mbox{Tr}_H\rangle \right )^a$ for some
exponent $a$.  Now $\euO_k[C_{2^n}]/\langle \mbox{Tr}_H\rangle$ is a
local ring with maximal ideal $\langle \sigma -1\rangle$ dividing
$2$. Therefore by Nakayama's Lemma any collection of elements in
$\mathcal{A}$ that serves as a $\euO_k/2\euO_k$-basis for
$\mathcal{C}/(\sigma-1)\mathcal{C}$ will serve as an
$\euO_k[C_{2^n}]/\langle \mbox{Tr}_H\rangle$-basis for $\mathcal{C}$.
This leaves us to show that $\mathcal{B}$ can be extended to a
$\euO_k/2\euO_k$-basis for the vector space
$\mathcal{C}/(\sigma-1)\mathcal{C}=\mathcal{A}/(\mathcal{A}^{H} +
(\sigma -1)\mathcal{A})$.  But since $\mbox{Tr}_{H}\mathcal{B}$ is a
basis for $\mbox{Tr}_{H}\mathcal{A}/\left
((\sigma-1)\mbox{Tr}_{H}\mathcal{A}+2\mathcal{A}^{H}\right )$, the
elements of $\mathcal{B}$ are linearly independent in
$\mathcal{A}/(\mathcal{A}^{H} + (\sigma -1)\mathcal{A})$ and therefore
span a subspace.
\end{proof}

\subsubsection{First Ramification Group of Order Four}
Let $|G_1|=4$.  This case is important because it illustrates the
utility of Lemma 2.2. (Recall that \S 2.2.1 and \S 2.2.3 are included
in this paper to motivate considerations in \S 3.)

\noindent {\em Step 1: Collect $|G_1|$ elements whose valuations are a
complete set of residues modulo $|G_1|$.}  We begin with the elements
used to determine the structure of ideals in $K_{n-1}$ (from \S 2.2.1),
namely $\alpha_m$ and $(\tilde{\sigma}+1)\alpha_m\in K_{n-1}$
(replacing $n$ by $n-1$, expressing $\sigma^{2^{n-2}}$ as
$\tilde{\sigma}$).  Note that the first ramification number of
$K_n/K_{n-2}$ is the (only) ramification number of $K_{n-1}/K_{n-2}$
(use \cite[pg 64 Cor]{serre:local} or switch to upper ramification
numbers \cite[IV \S 3]{serre:local}). So $v_n\left (\alpha_m \right
)=2v_{n-1}\left (\alpha_m \right )=2b_1+4m$ and $v_n((\tilde{\sigma}+1
)\alpha_m)=4b_1+4m$.  We have two elements of even valuation.  To get
elements with odd valuation, we apply Lemma 2.1(2). For each $X\in
K_{n-1}$, Lemma 2.1(2) gives us a preimage $\overline{X}\in K_n$ (under
the trace $\mbox{Tr}_{n,n-1}$), a preimage that satisfies
$v_n(\overline{X})=2v_{n-1}(X)-b_2$. So
$\mbox{Tr}_{n,n-1}\overline{X}=(\tilde{\sigma}^2+1)\overline{X}=X$.
The integers $v_n\left (\alpha_m \right )$, $v_n((\tilde{\sigma}+1
)\alpha_m)$, $v_n (\overline{\alpha_m} )=2b_1-b_2+4m$, $v_n
(\overline{(\tilde{\sigma}+1)\alpha_m})=4b_1-b_2+4m$ are a complete
set of residues modulo $4$.
\vspace*{3mm}

\noindent{\em Step 2: Collect elements with valuation in
$\{i,i+1,\ldots, i+v_n(2)-1\}$.}  To organize this process we use
Wyman's catalog of ramification numbers \cite{wyman}.  If $b_1\geq
e_0$, the second ramification number is uniquely determined,
$b_2=b_1+2e_0$.  If $b_1<e_0$, then either $b_2=3b_1$,
$b_2=4e_0-b_1$, or $b_2=b_1+4t$ for some $t$ with $b_1<2t<2e_0-b_1$
\cite[Thm 32]{wyman}.  In any case, we have the bound,
\begin{equation}
2b_1<b_2.
\end{equation}
Now for a given $m$, list the infinitely many elements,
$\alpha_{m+ke_0}$, $(\tilde{\sigma} +1)\alpha_{m+ke_0}$,
$\overline{\alpha_{m+ke_0}}$, $\overline{(\tilde{\sigma}
+1)\alpha_{m+ke_0}}$, in terms of increasing valuation. Replace
$\alpha_{m+ke_0}$ by $2^k\alpha_{m}$ and drop the subscripts $m$.  So for
$b_2> 4e_0-2b_1$, beginning at $\overline{\alpha}$, we have:
$$\cdots \longrightarrow\overline{\alpha}\overset{1}{\longrightarrow}
1/2(\tilde{\sigma}+1)\alpha \overset{2}{\longrightarrow}
\overline{(\tilde{\sigma} +1)\alpha} \overset{3}{\longrightarrow}
\alpha\overset{2}{\longrightarrow} \overline{2\alpha}\longrightarrow
\cdots $$ Each increase in valuation, denoted by
$\overset{x}{\longrightarrow}$, is justified as follows: For $x=1$,
the justification depends upon the case either $b_2> 4e_0-2b_1$ or
$b_2< 4e_0-2b_1$. For $x=2$, it is $b_2<4e_0$. For $x=3$, it is (2.2).
If $b_2< 4e_0-2b_1$, the list is as follows:
$$\cdots \longrightarrow\overline{\alpha}\overset{4}{\longrightarrow}
\overline{(\tilde{\sigma} +1)\alpha} \overset{3}{\longrightarrow}
\alpha \overset{4}{\longrightarrow}
(\tilde{\sigma}+1)\alpha\overset{1}{\longrightarrow}
\overline{2\alpha}\longrightarrow\cdots $$ Note $x=4$ is justified by
$b_1>0$.

Now collect those elements with valuation in $\{i,\ldots,
i+v_n(2)-1\}$.  This will provide us with an $\euO_{T(\zeta)}$-basis
for $\euP_n^i$. Begin with the smallest $m$ such that $i\leq
v_n(\alpha_m)$. Note then that $v_n (\overline{2(\tilde{\sigma}
+1)\alpha_m})<i+v_n(2)$.  Associated with this particular $m$ are four
elements in $\{i,\ldots, i+v_n(2)-1\}$.  They are listed in the first
row of the table below. Consider this interval to be a `window'.  As
we increase $m$, new elements appear ({\em e.g.}  $2X$) -- appearance
coincides with disappearance (namely of $X$).  Four elements are in
`view' always. There are four `views' (four sets).  We list the
`views' as rows under the appropriate heading. 
\vspace*{3mm}

\noindent $\mathcal{D}$: The $\euO_{T(\zeta)}$-basis
for $\euP_n^i$.
\begin{align}
{A}:\quad b_2< 4e_0-2b_1\quad\quad\quad &  \quad\quad\quad
{B}:\quad b_2> 4e_0-2b_1 \notag \\
\alpha\quad (\tilde{\sigma} +1)\alpha \quad
\overline{2\alpha}\quad\overline{2(\tilde{\sigma} +1)\alpha}\quad & \quad
\alpha \quad \overline{2\alpha}\quad (\tilde{\sigma}
+1)\alpha\quad\overline{2(\tilde{\sigma} +1)\alpha} \tag{$1$} \\
\overline{(\tilde{\sigma} +1)\alpha} \quad\alpha\quad (\tilde{\sigma} +1)\alpha \quad
\overline{2\alpha}\quad &\quad
\overline{(\tilde{\sigma} +1)\alpha}\quad \alpha \quad \overline{2\alpha}\quad
(\tilde{\sigma} +1)\alpha \tag{$2$} \\ 
\overline{\alpha}\quad\overline{(\tilde{\sigma} +1)\alpha} \quad\alpha\quad
(\tilde{\sigma} +1)\alpha\quad &\quad
1/2(\tilde{\sigma}
+1)\alpha\quad\overline{(\tilde{\sigma} +1)\alpha}\quad \alpha \quad
\overline{2\alpha} \tag{$3$} \\ 
1/2(\tilde{\sigma} +1)\alpha \quad
\overline{\alpha}\quad\overline{(\tilde{\sigma} +1)\alpha} \quad\alpha \quad &\quad
\overline{\alpha}\quad 1/2(\tilde{\sigma}
+1)\alpha\quad\overline{(\tilde{\sigma} +1)\alpha}\quad \alpha \tag{$4$}
\end{align}
Should we need to determine the subscripts (associated with a
particular `view'), we can easily do so: For example the four elements
listed in $A(1)$ and $B(1)$, appear for $m$ with $i\leq v_n(\alpha_m)$
and $v_n (\overline{2(\tilde{\sigma} +1)\alpha_m})\leq i+4e_0-1$. In
other words, $\lceil (i-2b_1)/4\rceil \leq m\leq \lceil
(i+b_2)/4\rceil -b_1-1$. 
\vspace*{3mm}

\noindent{\em Step 3: Identify a basis for the quotient module
$\euP_n^i/\euP_{n-1}^{\lceil i/2\rceil}$, and determine the precise
image of each basis element under the trace $\mbox{Tr}_{n,n-1}$ (in
terms of the basis for $\euP_{n-1}^{\lceil i/2\rceil}$).}  Observe
that $\euP_n^i/\euP_{n-1}^{\lceil i/2\rceil}$ is, in a natural way,
free over the principal ideal domain $\euO_{T(\zeta)}[G]/\langle
\tilde{\sigma}^2+1\rangle$.  We begin by identifying those elements
listed in $\mathcal{D}$, the $\euO_{T(\zeta)}$-basis from {\em Step
2}, that can serve as a $\euO_{T(\zeta)}[G]/\langle
\tilde{\sigma}^2+1\rangle$-basis.  Take $\mathcal{D}$ and partition it
into two sets. Let $\overline{\mathcal{D}}$ contain those elements
$\overline{X}$ with a bar. Let $\mathcal{D}_0$ contain those elements
$X$ without a bar.  So $\overline{\mathcal{D}}$ is an
$\euO_{T(\zeta)}$-basis for $\euP_n^i/\euP_{n-1}^{\lceil i/2\rceil}$,
and $\mathcal{D}_0$ is an $\euO_{T(\zeta)}$-basis for
$\euP_{n-1}^{\lceil i/2\rceil}$.  If we knew which elements from
$\overline{\mathcal{D}}$ provide us with $\euO_{T(\zeta)}[G]/\langle
\tilde{\sigma}^2+1\rangle$-basis for $\euP_n^i/\euP_{n-1}^{\lceil
i/2\rceil}$ we would be done, as it is easy to express the image
(under the trace $\mbox{Tr}_{n,n-1}$) of each element in
$\overline{\mathcal{D}}$ in terms of elements of $\mathcal{D}_0$
(there is a one--to--one correspondence).

Before we proceed further, note the following. We may assume without
loss of generality that for $\overline{X}\in\overline{\mathcal{D}}$,
$\mbox{Tr}_{n,n-1}\overline{X}\neq 0$ {\em if and only if}
$\overline{X}$ appears together with $X$ (for the same subscript $m$)
in $\mathcal{D}$. Clearly if $\overline{X}$ and $X$ appear together,
then $\mbox{Tr}_{n,n-1}\overline{X}=X\neq 0$. However when
$\overline{2X}$ and $X$ appear together, after a change of basis, we
may assume that $\mbox{Tr}_{n,n-1}\overline{2X}=0$.  The reason for
this is as follows: We can change an element of
$\overline{\mathcal{D}}$ by adding an element from $\mathcal{D}_0$ and
still have a $\euO_{T(\zeta)}$-basis for $\euP_n^i/\euP_{n-1}^{\lceil
i/2\rceil}$.  So whenever $\overline{2X}$ and $X$ appear together,
replace $\overline{2X}$ with $\overline{2X}-X$. Note
$\mbox{Tr}_{n,n-1}(\overline{2X}-X)=0$.  If we perform this change
throughout our basis, but relabel $\overline{2X}-X$ as
$\overline{2X}$, then we may continue to use the lists, $A(1)$--$A(4)$
and $B(1)$--$B(4)$, but assume that $\mbox{Tr}_{n,n-1}\overline{2X}=0$
if $\overline{2X}$ appears together with $X$.

Our next step will be to provide an $\euO_{T(\zeta)}[G]/\langle
\tilde{\sigma}^2+1\rangle$-basis for $\euP_n^i/\euP_{n-1}^{\lceil
i/2\rceil}$.  Consider those rows with an $\overline{X}$ such that
$\mbox{Tr}_{n,n-1}\overline{X}\neq 0$ (namely $A(2)$, $A(3)$, $A(4)$,
$B(2)$, $B(4)$). Let $\mathcal{S}\subseteq\overline{\mathcal{D}}$
denote the set of {\em left--most} $\overline{X}$ associated with
those rows.  So, for example, if $b_2+2b_1<4e_0$, then $\mathcal{S}$
is made up of the $\overline{(\tilde{\sigma} +1)\alpha_m}$ from
$A(2)$, and the $\overline{\alpha_m}$ from $A(3)$ and $A(4)$.  Verify
that $\mbox{Tr}_{n,n-1}\mathcal{S}$ is a
$\euO_{T(\zeta)}/2\euO_{T(\zeta)}$-basis for
$\mbox{Tr}_{n,n-1}\euP_n^i/ (
(\tilde{\sigma}-1)\mbox{Tr}_{n,n-1}\euP_n^i+ 2\euP_{n-1}^{\lceil
i/2\rceil} )$ (observe that $\mbox{Tr}_{n,n-1}\mathcal{S}$ generates
$\mbox{Tr}_{n,n-1}\euP_n^i/2\euP_{n-1}^{\lceil i/2\rceil}$ over
$\euO_{T(\zeta)}/2\euO_{T(\zeta)}[G]$).  Now use Lemma 2.2 to extend
$\mathcal{S}$ to $\mathcal{S}'$, an $\euO_{T(\zeta)}[G]/\langle
\tilde{\sigma}^2+1\rangle$-basis for $\euP_n^i/\euP_{n-1}^{\lceil
i/2\rceil}$.  Since $\euP_n^i/\euP_{n-1}^{\lceil i/2\rceil}$ has rank
$e_0$ over $\euO_{T(\zeta)}[G]/\langle \tilde{\sigma}^2+1\rangle$, we have
$|\mathcal{S}'|=e_0$.

This $\euO_{T(\zeta)}[G]/\langle \tilde{\sigma}^2+1\rangle$-basis,
$\mathcal{S}'$, possesses two important properties. First, it contains
$\mathcal{S}$. Second, without loss of generality we may assume that
the elements in $\mathcal{S}'-\mathcal{S}$ are killed by the trace
$\mbox{Tr}_{n,n-1}$. These two properties are shared with another set:
The set of {\bf all} {\em left--most} $\overline{X}$ (an $\overline{X}$ for
every value of $m$). Clearly the set of all {\em left--most}
$\overline{X}$ contains $\mathcal{S}$. Moreover, by an earlier
assumption, the compliment of $\mathcal{S}$ in the set of all {\em
left--most} $\overline{X}$ is mapped to zero under the trace. And so,
because the sets have the same cardinality (namely $e_0$), we can
identify them. Without loss of generality, assume that
$\mathcal{S}'$ {\em is} the set of all {\em left--most}
$\overline{X}$. This allows us to use the lists, $A(1)$--$A(4)$ and
$B(1)$--$B(4)$, in the `book-keeping' necessary for determining the
Galois module structure below.

At this point, we know that $\euP_n^i/\euP_{n-1}^{\lceil i/2\rceil}$
is free over $\euO_{T(\zeta)}[G]/\langle
\tilde{\sigma}^2+1\rangle$. Indeed, $\mathcal{S}'$ (the set of all
{\em left--most} $\overline{X}$) provides us a
$\euO_{T(\zeta)}[G]/\langle \tilde{\sigma}^2+1\rangle$-basis for
$\euP_n^i/\euP_{n-1}^{\lceil i/2\rceil}$. Of course, the
$\euO_{T(\zeta)}[G]$-structure of $\euP_{n-1}^{\lceil i/2\rceil}$ is
known from \S 2.2.1 (and can be read off of $\mathcal{D}_0$).  So a
description of the image of $\mathcal{S}'$ under $\tilde{\sigma}^2+1$
in terms of $\mathcal{D}_0$ will determine the Galois module
structure.  See \cite[\S 8]{curt}. {\em The Result}: For each $m$
associated with $A(1)$ or $B(1)$ we decompose off an
$\euO_{T(\zeta)}[G_1]$-summand of
$\euO_{T(\zeta)}\otimes_{\mathbb{Z}_2}\mathcal{I}$, for $A(2)$ or $B(2)$ we
get an $\euO_{T(\zeta)}\otimes_{\mathbb{Z}_2}\mathcal{H}$, for $A(3)$ we find
the group ring, $\euO_{T(\zeta)}[G_1]\cong
\euO_{T(\zeta)}\otimes_{\mathbb{Z}_2}\mathcal{G}$. But, for $B(3)$ we
decompose off the maximal order of $\euO_{T(\zeta)}[G_1]$,
$\euO_{T(\zeta)}\otimes_{\mathbb{Z}_2}\mathcal{M}$. For $A(4)$ and $B(4)$
there is $\euO_{T(\zeta)}\otimes_{\mathbb{Z}_2}\mathcal{L}$. All this and
counting determines the $\euO_{T(\zeta)}[G_1]$-module structure of
$\euP_n^i$ from which the $\euO_k[G_1]$-module structure can be
inferred. It also determines the module $\mathcal{M}_2(i,b_1,b_2)$.
To determine the $\euO_T[G]$-module structure (from which the
$\euO_k[G]$-module structure can be inferred), we need to take our
$\euO_{T(\zeta)}$-bases for $\euP_n^i$ and create $\euO_T$-bases.

\subsection{Partially Ramified Extensions}
Let $T_i$ denote the maximal unramified extension of $\mathbb{Q}_2$ contained
in $K_i$. So $T(\zeta)$ of the previous section can be expressed at
$T_n$, while $T=T_0$.  Recall the steps in \S 2.2.1.  We first
determined a $\euO_{T_n}$-basis $\mathcal{B}$ for $\euP_n^i$, one upon
which the action of $G_1$ could be explicitly followed.  Then noting
that we can identify $G/G_1$ with the Galois group for $T_n/T_0$, we
extended $\mathcal{B}$ to an $\euO_{T_0}$-basis for $\euP_n^i$. This
time the action of every element in the Galois group $G$ could be
followed.  What were the important ingredients in this process? It was
important that the elements of $\mathcal{B}$ lay in one--to--one
correspondence, via valuation, with the integers $i, \ldots,
i+v_n(2)-1$. Using this fact and the fact that for each $t$,
$\euP_n^t/\euP_n^{t+1}$ had a normal field basis over
$\euO_{T_0}/\euP_{T_0}$, we were able to make an $\euO_T$-basis for
$\euP_n^i$, namely $\mathcal{B}'=\cup_{\sigma^i\in
G/G_1}\sigma^i\mathcal{B}$.  At that point we were done. The
$\euO_T[G]$-structure could simply be read off of this basis. This is
not the case when $|G_1|=4$. Nor is it the case when $|G_1|=8$.  At
this point we still need to change our basis and use Nakayama's Lemma,
if only to determine $\euO_T[G_1]$-structure. We leave it to the
reader to check that this process of basis change `commutes' with the
process of extending our $\euO_{T_n}$-basis to an $\euO_{T_0}$-basis.
Simply follow the argument using elements of the form
$\sigma^t\alpha_m$, $\sigma^t(\tilde{\sigma}+1)\alpha_m,\ldots$ with
$t=0, \ldots 2^{[G:G_1]-1}$ instead of elements of the form
$\alpha_m$, $(\tilde{\sigma}+1)\alpha_m,\ldots$.  As a consequence,
the problem of determining the $\euO_T[G]$-module structure reduces to
the problem of determining the $\euO_{T_n}[G_1]$-module structure.

\section{Fully Ramified Cyclic Extensions of Degree Eight}
We consider fully ramified extensions $K_3/K_0$ with odd ramification
numbers.
\subsection{Outline}
Our discussion here is focused on the {\em unstably ramified} case,
$b_1<e_0$. (The {\em stably ramified} case will be addressed
separately in \S 3.4.)  Recall {\em Step 1} of \S 2.2.3 (in reference
to $K_2/K_0$).  But first note that the first two ramification numbers of
$K_3/K_0$ are the (only) two ramification numbers of $K_2/K_0$
\cite[pg 64 Cor]{serre:local}.  We
began with two elements, namely $\alpha, (\sigma+1)\alpha$ in the
subfield $K_1$. (The Galois relationship between them was explicit.)
Then we created $\overline{\alpha}, \overline{(\sigma+1)\alpha}\in
K_2$, preimages under the trace $\mbox{Tr}_{2,1}$.  In this section,
we will start with these four elements from $K_2$ and use
Lemma 2.1(2) to find further preimages: of $\alpha, (\sigma+1)\alpha,
\overline{\alpha}, \overline{(\sigma+1)\alpha}$ under
$\mbox{Tr}_{3,2}$.  To avoid confusion (confusion resulting from
additional bars denoting a preimage under $\mbox{Tr}_{3,2}$), we
relabel. Let $\alpha:=\overline{\alpha}$ and let
$\rho:=\overline{(\sigma+1)\alpha}$.  So the four elements in $K_2$
are labeled $\alpha, (\sigma^2+1)\alpha, \rho,
(\sigma+1)(\sigma^2+1)\alpha$ (instead of $\overline{\alpha}, \alpha,
\overline{(\sigma+1)\alpha}, (\sigma+1)\alpha$ respectively).  The
eight resulting elements (four from $K_2$ along with their preimages)
lie in one--to--one correspondence with the residues modulo $8$.

We would have accomplished all that was accomplished in {\em Step 1}
from \S 2.2.3 if we knew the Galois relationships among $\alpha,
(\sigma^2+1)\alpha, \rho, (\sigma+1)(\sigma^2+1)\alpha$ explicitly.
We need an explicit relationship between $\alpha$ and $\rho$. This is
accomplished in \S 3.2.2 through a list of technical results --
generalizations of Lemma 2.1. Note that $\rho$ is an `approximation'
to $(\sigma+1)\alpha$ -- they have the same image under the trace
$\mbox{Tr}_{2,1}$.  Our results describe their difference, the `error'
in this `approximation'.

As a prerequisite for the technical results of \S 3.2.2, and in
preparation for the analog of {\em Step 2} from \S 2.2.3 we use a
result of Fontaine to provide a catalog of ramification numbers in \S
3.2.1. We are then ready for {\em Step 2}: First we order the eight
elements (that we inherit from {\em Step 1}) in terms of increasing
valuation.  This is accomplished in \S 3.3.  There are eight orderings
-- eight cases. The result is eight different bases, listed as $A$ --
$H$ (as opposed to just two in $\mathcal{D}$ from \S 2.2.3). For the
convenience of the reader, they are listed in Appendix B.

We are now ready for the analog of {\em Step 3} from \S 2.2.3.  We are
ready to determine those elements in each $\euO_T$-basis that serve as
an $\euO_T[G]/\langle\mbox{Tr}_{3,2}\rangle$-basis, $\mathcal{S}$, for
$\euP_3^i/\euP_2^{\lceil i/2\rceil}$.  We will then be able to
describe the image, $\mbox{Tr}_{3,2}S$, in terms of our $\euO_T$-basis
for $\euP_2^{\lceil i/2\rceil}$ (or more to the point, explicitly in
terms of $\euO_T[G]$-generators for $\euP_2^{\lceil i/2\rceil}$). To
do all this we will need, as in \S 2.2.3, to perform certain basis
changes. The processes are similar, but there are a few very important
differences.  For the convenience of the reader, the results of this
step are summarized in \S 3.3.1. The steps are then spelled out in \S
3.3.2 -- \S 3.3.5.  The structure of $\euM_3(i,b_1, b_2,b_3)$ (given
in Tables 1 and 2) can then be read off of the bases in Appendix
B. Note however, that we still need to determine the structure under
$b_1\geq e_0$ (part of Case $A$).  This situation is addressed in \S
3.3.4.

\subsection{Preliminary Results} We catalog the ramification triples and
generalize Lemma 2.1, describing the difference $\rho-(\sigma+1)\alpha$.
\subsubsection{Ramification Triples}
There is {\em stability} and {\em instability}.

\begin{thm}[{\cite[Prop 4.3]{fontaine}}] $\mbox{}$
\newline Stability:
$$b_1\geq e_0 \Rightarrow b_2=b_1+2e_0,\qquad\mbox{and}\qquad b_1+b_2\geq 2e_0
\Rightarrow b_3=b_2+4e_0.$$ 
Instability:
$$b_1< e_0 \Rightarrow 3b_1\leq b_2\leq 4e_0-b_1,\qquad b_1+b_2< 2e_0
\Rightarrow 3b_2+2b_1\leq b_3\leq 8e_0-b_2-2b_1.$$ In particular, when
$b_1<e_0$, either $b_2=3b_1$, $b_2=4e_0-b_1$, or $b_2=b_1+4t$ for
$b_1<2t<2e_0-b_1$, while if $b_1+b_2<2e_0$, then either
$b_2=3b_2+2b_1$, $b_2=8e_0-b_2-2b_1$, or $b_3=8s-b_2+2b_1$ for
$b_2<2s<2e_0-b_1$.
\end{thm}

Plot these ramification triples $(b_1, b_2, b_3)$ in $\Re^3$, and
project this plot to the first two coordinates, $(x, y, z)\rightarrow
(x, y, 0)$, thus creating Figure 1 (next page). This projection is
partly a line: for $b_1\geq e_0$, each point $(b_1, b_2)$ is
restricted to $b_2=b_1+2e_0$. It is partly a triangular region: for
$b_1<e_0$, each point $(b_1, b_2)$ is bound between the lines
$b_2=3b_1$ and $b_2=4e_0-b_1$.  The significance of the regions $A, B,
C,\ldots $ will be explained later. Note that for points, $(b_1,
b_2)$, above the line $b_2=-b_1+2e_0$, the plot of the $(b_1, b_2,
b_3)$ in $\Re^3$ will be a plane -- $b_3$ is uniquely determined.

In Figure 2 we have plotted a slice, at a particular value of $b_1$,
through our plot of ramification triples in $\Re^3$.  Part of this
slice is a line -- when $b_3$ is uniquely determined.  Thus the line
from $(2e_0-b_1, 6e_0-b_1)$ to $(4e_0-b_1, 8e_0-b_1)$.  Indeed, as
drawn, Figure 2 implicitly assumes that the slice was taken at $b_1$
for $b_1<e_0/2$. Otherwise there would be no triangular region.
Observe that in Figure 1, the lines $b_2=2e_0-b_1$ and $b_2=3b_1$
intersect at $b_1=e_0/2$. If $b_1\geq e_0/2$, the third ramification
number is uniquely determined by $b_2$.  The triangular region bound
by the lines $b_2=3b_1$, $b_3=3b_2+2b_1$ and $b_3=8e_0-b_2-2b_1$
exists only for $b_1<e_0/2$.

Because the ramification numbers are odd, the triangular part of
Figure 1 can be partitioned as follows:
\begin{itemize}
\item[$A$.] $4e_0-4b_1/3<b_2$
\item[$B$.] $4e_0-2b_1< b_2 < 4e_0-4b_1/3$
\item[$C$.] $4e_0-4b_1< b_2 < 4e_0-2b_1$ and $b_2>(4e_0+4b_1)/3$
\item[$D$.] $4e_0-4b_1< b_2 < 4e_0-2b_1$ and $b_2<(4e_0+4b_1)/3$
\item[$E$.] $2e_0-b_1\leq b_2 < 4e_0-4b_1$ and $b_2>(4e_0+4b_1)/3$
\item[$F$.] $2e_0-b_1\leq b_2 < 4e_0-4b_1$ and $b_2<(4e_0+4b_1)/3$
\end{itemize}

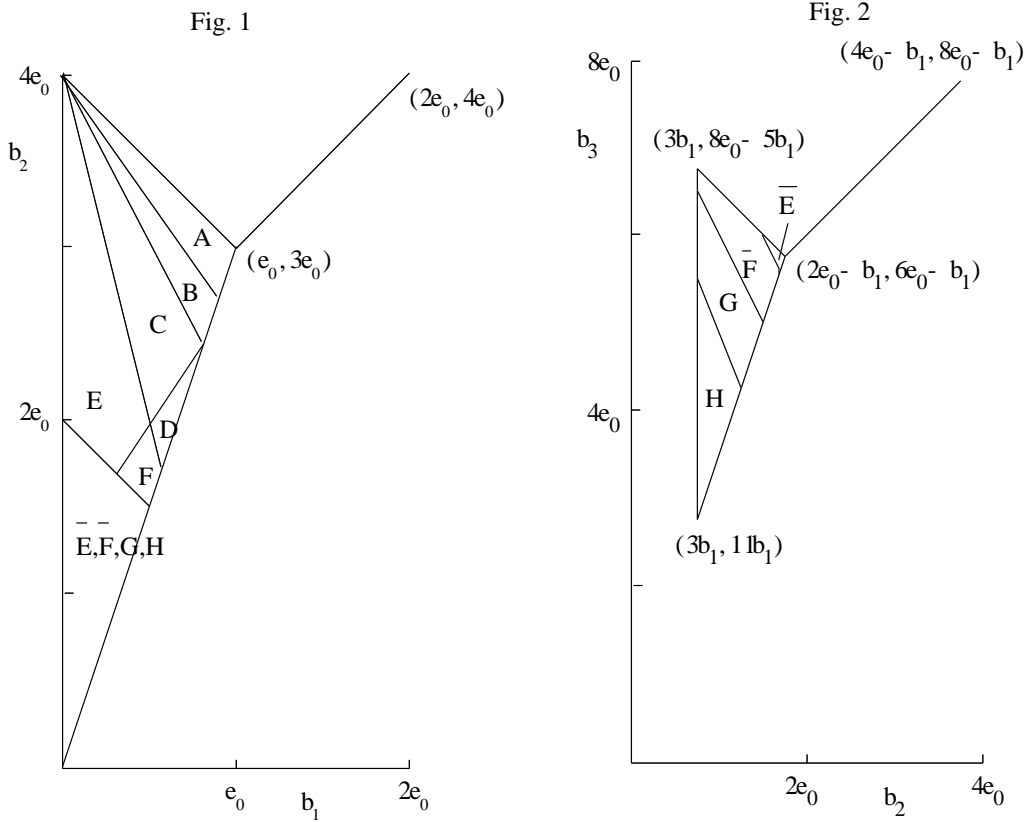
\vspace*{6mm}

\hspace*{1.25in}\hspace{-4cm}\begin{picture}(0,0)%
\includegraphics{picture}%
\end{picture}%
\setlength{\unitlength}{3947sp}%
\begingroup\makeatletter\ifx\SetFigFont\undefined
\def\x#1#2#3#4#5#6#7\relax{\def\x{#1#2#3#4#5#6}}%
\expandafter\x\fmtname xxxxxx\relax \def\y{splain}%
\ifx\x\y   
\gdef\SetFigFont#1#2#3{%
  \ifnum #1<17\tiny\else \ifnum #1<20\small\else
  \ifnum #1<24\normalsize\else \ifnum #1<29\large\else
  \ifnum #1<34\Large\else \ifnum #1<41\LARGE\else
     \huge\fi\fi\fi\fi\fi\fi
  \csname #3\endcsname}%
\else
\gdef\SetFigFont#1#2#3{\begingroup
  \count@#1\relax \ifnum 25<\count@\count@25\fi
  \def\x{\endgroup\@setsize\SetFigFont{#2pt}}%
  \expandafter\x
    \csname \romannumeral\the\count@ pt\expandafter\endcsname
    \csname @\romannumeral\the\count@ pt\endcsname
  \csname #3\endcsname}%
\fi
\fi\endgroup
\begin{picture}(6326,5166)(1289,-4867)
\end{picture}

Assuming that $b_1<e_0/2$, there is a triangular
region in Figure 2. This can be partitioned into the following cases:
\begin{itemize}
\item[$\overline{E}$.] $8e_0+4b_1-2b_2 < b_3$
\item[$\overline{F}$.] $8e_0-2b_2<b_3<8e_0+4b_1-2b_2$
\item[$G$.] $8e_0-4b_1-2b_2 < b_3 < 8e_0-2b_2$
\item[$H$.] $b_3<8e_0-4b_1-2b_2$
\end{itemize}

Note that if $b_1>8e_0/17$, region $G$ is empty; if $b_1>8e_0/21$,
region $H$ is empty; if $b_1>8e_0/28$, region $\overline{E}$ is empty.
So as drawn, we have assumed that $b_1<2e_0/7$.  If however the slice
were taken for a value $8e_0/17<b_1<8e_0/16$, note that the triangular
region would consist of only one case, namely $\overline{F}$. The
relationship between $E$, $F$ and $\overline{E}$, $\overline{F}$ will
be explained in \S 3.3.

\subsubsection{Technical Lemmas}
The difference $\rho-(\sigma+1)\alpha$ depends upon ramification.
\paragraph{\em Unstable Ramification}
Assume that $b_1<e_0$.  These results may be thought of as
consequences of indirect `routes' from $\alpha$ to $\rho$. For example,
we may begin with $\alpha\in K_2$, create
$(\sigma^2+1)\alpha$, then $(\sigma+1)(\sigma^2+1)\alpha$ and let
$\rho$ be the inverse image of $(\sigma+1)(\sigma^2+1)\alpha$ under
$\mbox{Tr}_{2,1}$. This results in an expression for the
$\rho-(\sigma+1)\alpha$.

\begin{lem} 
If $b_2\equiv b_1\bmod 4$ {\rm (}equivalently $3b_1<b_2<4e_0-b_1${\rm )}, let $t=(b_2-b_1)/4$.
There are elements $\alpha_m\in K_2$ with
$v_2(\alpha_m)=b_2+4m$, such that
$$\rho_m=(\sigma +1)\alpha_m+ (\sigma^2\pm 1)\alpha_{m-t}$$ has
valuation $v_2(\rho_m)=b_2+2b_1+4m$. The `$+$' or `$-$' depends on
our needs.
\end{lem}
\begin{proof}
Let $\alpha\in K_2$ with valuation,
$v_2(\alpha)\equiv b_2\bmod 4$.  Using Lemma 2.1,
$v_2((\sigma+1)\alpha)=v_2(\alpha)+b_1$,
$v_2((\sigma^2+1)\alpha)=v_2(\alpha)+b_2$. Since
$(\sigma^2+1)\alpha\in K_1$ and
$v_1((\sigma^2+1)\alpha)=(v_2(\alpha)+b_2)/2\equiv b_2\bmod 2$,
$v_1((\sigma+1)(\sigma^2 +1)\alpha)=(v_2(\alpha)+b_2)/2+b_1$.  Using
Lemma 2.1(2), there is a $\rho\in K_2$ with
$v_2(\rho)=v_2(\alpha)+2b_1$ such that
$(\sigma^2+1)\rho=(\sigma+1)(\sigma^2 +1)\alpha$.  Since
$(\sigma^2+1)\left [\rho-(\sigma+1)\alpha\right ]=0$. Using Lem
2.1(3), there is a $\theta\in K_2$ with
$v_2(\theta)=(v_2(\alpha)-b_2)+b_1$ and
$\rho=(\sigma+1)\alpha+(\sigma^2-1)\theta$.  Since $b_1 < e_0$,
$v_2(2\theta)>v_2(\rho)$. We may replace $\rho$ by
$\rho':=\rho+2\theta$ (they have the same valuation), and get
$\rho'=(\sigma+1)\alpha+(\sigma^2+1)\theta$.  Once
$\alpha_m$ is chosen, we let $\alpha_{m-t}:=\theta$.
\end{proof}

\begin{lem} If $b_2\equiv -b_1\bmod 4$ {\rm (}equivalently $b_2=3b_1$ or
$b_2=4e_0-b_1${\rm )}, let $s:=(b_2+b_1)/4$.
There are elements $\alpha_m\in K_2$
with $v_2(\alpha_m)=b_2+4m$,
such that
$$\rho_m=(\sigma +1)\alpha_m + (\sigma+1)(\sigma^2+1)\alpha_{m-s}$$
has valuation, $v_2(\rho_m)=2b_2-b_1+4m$. 
Note if $b_2=3b_1$, $v_2(\rho_m)=b_2+2b_1+4m$.
\end{lem}
\begin{proof}
There is a $\tau\in K_0$ with
$v_0(\tau)=(b_2-b_1)/2$. Using Lemma 2.1(2), let
$\rho\in K_2$ with $v_2(\rho)=b_2-2b_1$ such that $(\sigma ^2+1)\rho
=\tau$.  Clearly $(\sigma^2+1)\cdot (\sigma -1)\rho=0$, so there is a
$\theta\in K_2$ with $v_2(\theta)=-b_1$ such that $(\sigma
-1)\rho=(\sigma^2-1)\theta$.  Since $(\sigma -1)\cdot [\rho - (\sigma
+1)\theta]=0$, $\tau':=\rho-(\sigma+1)\theta$ is a unit in $K_0$.  Let
$\rho'=\rho/\tau'$ and $\theta'=\theta/\tau'$, so $1=\rho'-(\sigma
+1)\theta'$.  Now let $\beta=(\sigma
+1)(\sigma^2+1)\theta'$. Clearly $(\sigma^2+1)\theta'\in K_1$ and
$v_1((\sigma^2+1)\theta')=(b_2-b_1)/2$ is odd. Therefore
$v_2(\beta)=b_2+b_1$.  Replacing $1$ with the expression, $(\sigma
+1)(\sigma^2+1)(\theta'/\beta)$, yields
\begin{equation}
\rho'=(\sigma +1)\theta'+(\sigma
+1)(\sigma^2+1)(\theta'/\beta).
\end{equation}
By choosing $\tau\in K_0$ with other valuations, the result follows.
\end{proof}

Unfortunately, if $b_2=4e_0-b_1$ then $s=e_0$ (valuation can not
distinguish between $\alpha_m/2$ and $\alpha_{m-s}$).  To avoid this
confusion, we include the following.
\begin{lem} 
Let $b_2=4e_0-b_1$.  There are $\alpha_m\in K_2$ with
$v_2(\alpha_m)=b_2+4m$, so
$$\rho_m:=(\sigma +1)\alpha_m - \frac{1}{2}(\sigma +1)(\sigma^2+1)\alpha_m
+ \frac{1}{2}(\sigma +1)(\sigma^2+1)\alpha_{m+e_0-b_1}$$
has valuation, $v_2(\rho_m)=2b_2-b_1+4m$.
\end{lem}
\begin{proof}
From (3.1) we have $\rho'=(\sigma +1)\theta'+(\sigma
+1)(\sigma^2+1)(\theta'/\beta)$.  Apply $(\sigma^2+1)/\beta$ to
both sides. So $(\sigma^2+1)(\rho'/\beta)=1+2/\beta$.  Since
$v_2((\sigma^2+1)\rho')=8e_0-4b_1$ and $v_2(\beta )=4e_0$, then
$v_0(1+2/\beta)=e_0-b_1$.  Replace $\theta'/\beta$ with
$(1/2)\cdot [-\theta'+\theta'(1+2/\beta)]$, and distribute $(\sigma
+1)(\sigma^2+1)$. 
\end{proof}
\begin{rem} 
Note $(\sigma-1)\rho_m=(\sigma^2 -1)\alpha_m$ and
$(\sigma^2+1)\rho_m=(\sigma +1)(\sigma^2+1)\alpha_{m+e_0-b_1}$, using
Lemma 3.4.  Apparently, $\rho_m$ is `torn' between $\alpha_m$ and
$\alpha_{m+e_0-b_1}$.  We chose to emphasize $\rho_m$'s tie to $\alpha_m$.  If we
relabel $\rho_{m-e_0+b_1}$ as $\rho_m$ (keep the $\alpha_m$ the
same), Lemma 3.4 reads
$$\rho_m:=(\sigma +1)\alpha_{m-e_0+b_1} - \frac{1}{2}(\sigma
+1)(\sigma^2+1)\alpha_{m-e_0+b_1} + \frac{1}{2}(\sigma
+1)(\sigma^2+1)\alpha_m$$ has valuation, $v_2(\rho_m)=b_2+2b_1+4m$ --
thus tying $\rho_m$ to $(1/2)(\sigma
+1)(\sigma^2+1)\alpha_m$.  This valuation of $\rho_m$ is as in
Lemmas 3.2 and 3.3 (for
$b_2=3b_1$).
\end{rem}

\paragraph{\em Stably Ramified Extensions} 
Assume that $b_1\geq e_0$.  The results may be seen as direct routes
from $\alpha$ to $\rho$. We create $\rho$ immediately from $(\sigma
+1)\alpha\in K_2$.  For discussion and generalization, see
\cite{elder:bord}.
\begin{lem} Let $b_1>e_0$.
For every odd integer, $a$, there are elements $\alpha, \rho\in K_2$
with $v_2(\alpha)=a$, $v_2(\rho)=a+(b_2-b_1)$. such that
$$(\sigma + 1)\alpha -\rho = \mu \in K_1,$$ with $v_2(\mu)=v_2(\alpha
)+b_1$.  Furthermore $\mu\in K_0$ for $v_2(\mu)=v_2(\alpha )+b_1\equiv
0\bmod 4$.
\end{lem}
\begin{proof}
Since $v_2((\sigma +1)\alpha)=v_2(\alpha)+b_1$ is even, we may express
$(\sigma +1)\alpha$ as a sum $\mu +\rho$ with $\mu\in K_1$, $\rho\in
K_2$, $v_2(\mu)=v_2(\alpha )+b_1$ and odd $v_2(\rho)$. Apply $(\sigma
-1)$. So $(\sigma ^2-1)\alpha = (\sigma -1)\mu + (\sigma -1)\rho$.
Since $b_2=b_1+2e_0<3b_1$, $v_2((\sigma ^2-1)\alpha)=v_2(\alpha )+b_2
<v_2(\alpha )+3b_1\leq v_2((\sigma -1)\mu)$. So $v_2((\sigma
^2-1)\alpha)=v_2((\sigma -1)\rho )$ and
$v_2(\rho)=v_2(\alpha)+(b_2-b_1)$.  If $v_2(\mu)\equiv 0 \bmod 4$, we
may choose $\alpha$ so that $\mu\in K_0$.  Pick a $\mu^*\in K_0$ with
$v_2(\mu^*)=v_2(\mu)$.  Relabel $\alpha$ as $\alpha_0$. Choose
$\alpha_i\in K_2$ with $v_2(\alpha_i)=v_2(\alpha_0)+2i$. As before,
generate $\mu_i$ and $\rho_i$ with $\alpha_i=\mu_i+\rho_i$. Clearly
$\mu^*=\sum_{i=0}^\infty a_i\mu_i$ for some units $a_i\in K_0$.  Let
$\alpha^*=\sum_{i=0}^\infty a_i\alpha_i$ and $\rho^*=\sum_{i=0}^\infty
a_i\rho_i$.
\end{proof}

\begin{lem} Let $b_1=e_0$ be odd.
For every odd integer, $a$, there are elements $\alpha, \rho\in K_2$
with $v_2(\alpha)=a$, $v_2(\rho)=a+(b_2-b_1)$ such that
\begin{align*}
(\sigma - 1)\alpha -\rho &= \mu_1\in K_1 \mbox{ if $a\equiv e_0\bmod 4$},\\
(\sigma +1)\alpha -\rho &= \mu_0\in K_0  \mbox{ if $a\equiv 3e_0\bmod 4$}.
\end{align*}
with
$v_2(\mu_i)=v_2(\alpha )+b_1$.
\end{lem}
\begin{proof}
Let $\tau\in K_0$ be a unit. From Lemma 2.1(2), there is a $\rho\in K_2$
with $v_2(\rho)=-b_2$ and $(\sigma^2+1)\rho=\tau$. So
$(\sigma^2+1)\cdot(\sigma -1)\rho=0$. Use Lemma 2.1(3) to find
$\theta\in K_2$ with $v_2(\theta)=b_1-2b_2$ and
$(\sigma^2-1)\theta=(\sigma -1)\rho$.  For $a\equiv e_0\bmod 4$, we
may assume that $\alpha=\rho\pi_0^m$ for some $m$.  Let
$\mu_1=(\sigma^2+1)\theta\pi_0^m\in K_1$ and $\rho= -2\theta\pi_0^m\in
K_2$. The statement follows.  For $a\equiv 3e_0\bmod 4$,
$(\sigma^2-1)\theta=(\sigma -1)\rho$ can be interpreted to mean that
$\rho - (\sigma +1)\theta\in K_0$.  Multiplying by an appropriate
power of $\pi_0$, we let $\alpha = \theta\pi_0^m$, $\mu_0=-(\rho -
(\sigma +1)\theta)\pi_0^m\in K_0$ and $\rho=\rho\pi_0^m\in K_2$.
\end{proof}

\subsection{The Galois module structure under {\em unstable} ramification}
Assume $b_1<e_0$.
First we determine the $\euO_T$-bases
in Appendix B. 
From Lemmas 3.2, 3.3, 3.4 we have
$\alpha_m$, $\rho_m$, $(\sigma
^2+1)\alpha_m$, $(\sigma+1)(\sigma ^2+1)\alpha_m \in K_2$,
 with valuations (measured in $v_2$) for every residue class modulo
$4$.  Recall $v_2(\alpha_m)=b_2+4m$,
$v_2((\sigma^2+1)\alpha_m)=2b_2+4m$, $v_2((\sigma+1)(\sigma
^2+1)\alpha_m)=2b_2+2b_1+4m$ and 
$v_2(\rho_m)=
8e_0-3b_1+4m$ if  $b_2= 4e_0-b_1$, otherwise
$v_2(\rho_m)=b_2+2b_1+4m$.
Using Lemma 2.1(2) we determine elements
$\overline{\alpha_m}$,
$\overline{\rho_m}$,
$\overline{(\sigma^2+1)\alpha_m}$,
$\overline{(\sigma+1)(\sigma^2+1)\alpha_m}\in K_3$, with
$(\sigma^4+1)\overline{X}=X$ and $v_3(\overline{X})=2v_2(X)-b_3$.
These eight elements have valuations (measured in $v_3$) in
one--to--one correspondence with the residue classes modulo $8$. By
varying $m$, it is possible to choose those with valuation $i\leq
v_3(x)<8e_0+i$.

To organize this process, we list these elements in
terms of increasing valuation.  There are eight orderings -- eight cases.
In each case $X$ (or $\overline{X}$), an increase in
valuation is denoted by an arrow, $\longrightarrow$, and justified by an
inequality assigned a number.  Numbers above an arrow apply to
$X$. Numbers below the arrow apply to $\overline{X}$.  As
we see below, the ordering of the elements in $\overline{E}$ is the
same as in $E$ (also in $\overline{F}$ as in $F$). This explains the use
of similar notation.

\begin{multline*}\tag*{$A$.}
\rho \overset{1}{\longrightarrow} \overline{2\rho}
\overset{2}{\longrightarrow} (\sigma^2+1)\alpha
\overset{1}{\longrightarrow} \overline{2(\sigma^2+1)\alpha}
\overset{1}{\longrightarrow} 2\alpha \overset{1}{\longrightarrow} \\
\overline{4\alpha} \overset{6}{\longrightarrow} (\sigma
+1)(\sigma^2+1)\alpha \overset{1}{\longrightarrow} \overline{2(\sigma
+1)(\sigma^2+1)\alpha} \overset{4}{\longrightarrow} 2\rho
\end{multline*}
In Case $A$, the valuation of $\rho_m$ depends upon whether or not
$b_2=4e_0-b_1$. If $b_2=4e_0-b_1$, $0<b_1$ justifies $2$ while
$b_1<2e_0$ justifies $4$. All other increases, including $2$ and $4$
for $b_2\neq 4e_0-b_1$, are justified by the inequalities listed
below.  

In Cases $B$ through $H$, there is only one valuation of $\rho_m$.
\begin{multline*}\tag*{$B$.}
\rho \overset{1}{\longrightarrow} \overline{2\rho}
 \overset{2}{\longrightarrow} (\sigma^2+1)\alpha
 \overset{1}{\longrightarrow} \overline{2(\sigma^2+1)\alpha}
 \overset{4}{\longrightarrow} 2\alpha \overset{5}{\longrightarrow} \\
 (\sigma +1)(\sigma^2+1)\alpha
 \overset{6'}{\longrightarrow}\overline{4\alpha}
 \overset{5}{\longrightarrow} \overline{2(\sigma
 +1)(\sigma^2+1)\alpha} \overset{4}{\longrightarrow} 2\rho
\end{multline*}
\begin{multline*}\tag*{$C$.}
\rho \overset{1}{\longrightarrow} \overline{2\rho}
\overset{2}{\longrightarrow} (\sigma^2+1)\alpha
\overset{1}{\longrightarrow} \overline{2(\sigma^2+1)\alpha}
\overset{7}{\longrightarrow} (\sigma +1)(\sigma^2+1)\alpha
\overset{5'}{\longrightarrow} \\ 2\alpha \overset{7}{\longrightarrow}
\overline{2(\sigma +1)(\sigma^2+1)\alpha}
\overset{5'}{\longrightarrow}\overline{4\alpha}\overset{7}{\longrightarrow}2\rho
\end{multline*}
\begin{multline*}\tag*{$D$.}
\rho\overset{9}{\longrightarrow} (\sigma^2+1)\alpha
\overset{2'}{\longrightarrow} \overline{2\rho}
\overset{9}{\longrightarrow} \overline{2(\sigma^2+1)\alpha}
\overset{7}{\longrightarrow} (\sigma +1)(\sigma^2+1)\alpha
\overset{5'}{\longrightarrow} \\ 2\alpha \overset{7}{\longrightarrow}
\overline{2(\sigma +1)(\sigma^2+1)\alpha}
\overset{5'}{\longrightarrow}\overline{4\alpha}\overset{7}{\longrightarrow}2\rho
\end{multline*}
\begin{multline*}\tag*{$E=\overline{E}$.}
\rho\underset{11}{\overset{7'}{\longrightarrow}}
\overline{2\alpha}\underset{8}{\overset{8}{\longrightarrow}}
\overline{2\rho} \underset{13}{\overset{2}{\longrightarrow}}
(\sigma^2+1)\alpha \underset{8}{\overset{8}{\longrightarrow}} (\sigma
+1)(\sigma^2+1)\alpha \underset{11}{\overset{7'}{\longrightarrow}} \\
\overline{2(\sigma^2+1)\alpha}
\underset{8}{\overset{8}{\longrightarrow}} \overline{2(\sigma
+1)(\sigma^2+1)\alpha} \underset{14}{\overset{7'}{\longrightarrow}}
2\alpha \underset{8}{\overset{8}{\longrightarrow}} 2\rho
\end{multline*}
\begin{multline*}\tag*{$F=\overline{F}$.}
\rho\underset{11}{\overset{7'}{\longrightarrow}}
\overline{2\alpha}\underset{12}{\overset{10}{\longrightarrow}}
(\sigma^2+1)\alpha \underset{13'}{\overset{2'}{\longrightarrow}}
\overline{2\rho} \underset{12}{\overset{10}{\longrightarrow}} (\sigma
+1)(\sigma^2+1)\alpha \underset{11}{\overset{7'}{\longrightarrow}} \\
\overline{2(\sigma^2+1)\alpha}
\underset{8}{\overset{8}{\longrightarrow}} \overline{2(\sigma
+1)(\sigma^2+1)\alpha} \underset{14}{\overset{7'}{\longrightarrow}}
2\alpha \underset{8}{\overset{8}{\longrightarrow}} 2\rho
\end{multline*}
\begin{multline*}\tag*{$G$.}
\rho
\overset{9}{\longrightarrow} (\sigma^2+1)\alpha 
\overset{12'}{\longrightarrow}\overline{2\alpha}
\overset{15}{\longrightarrow} (\sigma +1)(\sigma^2+1)\alpha 
\overset{12'}{\longrightarrow} \overline{2\rho}  
\overset{9}{\longrightarrow} \overline{2(\sigma^2+1)\alpha} \overset{8}{\longrightarrow} \\
 \overline{2(\sigma +1)(\sigma^2+1)\alpha}  \overset{14}{\longrightarrow}
2\alpha 
\overset{8}{\longrightarrow} 2\rho
\end{multline*}
\begin{multline*}\tag*{$H$.}
\rho
\overset{9}{\longrightarrow} (\sigma^2+1)\alpha 
\overset{8}{\longrightarrow} (\sigma +1)(\sigma^2+1)\alpha 
\overset{15'}{\longrightarrow}\overline{2\alpha}
\overset{8}{\longrightarrow} \overline{2\rho}  
\overset{9}{\longrightarrow} \overline{2(\sigma^2+1)\alpha} \overset{8}{\longrightarrow} \\
 \overline{2(\sigma +1)(\sigma^2+1)\alpha}  \overset{14}{\longrightarrow}
2\alpha 
\overset{8}{\longrightarrow} 2\rho
\end{multline*}
{\em Numbered Inequalities}:
$(1)$ $b_1<2e_0$, $b_2<4e_0$, $b_3<8e_0$.
$(2)$ $3b_2>4e_0+4b_1$. $(2')$ $3b_2<4e_0+4b_1$.
$(3)$ $4e_0-4b_1<3b_2$ (true for $A$--$F$ since $b_2\geq 2e_0-b_1$).
$(4)$ $2b_2<b_3$.
$(5)$ $4e_0-2b_1<b_2$. $(5')$ $4e_0-2b_1>b_2$.
$(6)$ $4e_0-4b_1/3<b_2$. $(6')$ $4e_0-4b_1/3>b_2$.
$(7)$ $4e_0-4b_1<b_2$. $(7')$ $4e_0-4b_1>b_2$.
$(8)$ $b_1>0$.
$(9)$ $b_2>2b_1$.
$(10)$ $b_2>4e_0/3$ (true for $A$--$F$, since $b_2\geq 3e_0/2$).
$(11)$ Since $b_2>2b_1$ and $b_3\leq 8e_0-2b_1-b_2$, $b_3<8e_0-4b_1$.
$(12)$ $8e_0-2b_2<b_3$. $(12')$ $8e_0-2b_2>b_3$.
$(13)$ $8e_0+4b_1-2b_2<b_3$. $(13')$ $8e_0+4b_1-2b_2>b_3$.
$(14)$ Since $b_2>2b_1$ and $b_3\geq 3b_2+2b_1$, $b_3>2b_1+4b_2$.
$(15)$ $8e_0-4b_1-2b_2<b_3$. $(15')$ $8e_0-4b_1-2b_2>b_3$.

We leave it to the reader to verify Appendix B.

\subsubsection{Summary: Results of Basis Changes and Nakayama's Lemma}

\paragraph{\em Basis Changes}
Except in four
rows,
\begin{equation}
C(2), D(2), E(2), F(2),
\end{equation}
we find we may change the $\euO_T$-bases in Appendix B so that the
Galois action upon each basis is {\bf\em as if} $\rho$ and
$\overline{\rho}$ had been everywhere replaced by $(\sigma +1)\alpha$
and $\overline{(\sigma +1)\alpha}$.  In the four exceptional cases
there are nontrivial Galois relationships among the basis
elements. This is explained in \S 3.3.5.

\paragraph{\em Nakayama's Lemma}
We find, without loss of generality, that the set $\mathcal{S}$ of
`left--most' elements $\overline{X}$ (as in $\mathcal{S}'$ of \S
2.2.3) from each basis in Appendix B will serve as a
$\euO_T[G]/\langle\mbox{Tr}_{3,2}\rangle$-basis for
$\euP_3^i/\euP_2^{\lceil i/2\rceil}$, {\em except} that $\mathcal{S}$
contains both $\overline{(\sigma+1)(\sigma^2+1)\alpha }$ and
$\overline{2\alpha}$ in $B(3)$, $C(3)$, $D(3)$.

At this point, the reader can skip the verification of these
assertions, ignore Cases $C$ through $F$, replace $\rho$ with $(\sigma
+1)\alpha$, and lift the Galois module structure off of the bases
listed in Appendix B.  See \cite[\S 8]{curt} The result of the readers
effort will be the statement of our main result in every case except
those associated with (3.2).

\subsubsection{Trivial Difference}
The elements $\alpha_m$, $\rho_m$ (or $\rho_m$, $2\alpha_m$)
from each basis in Appendix B provide a $\euO_T$-basis for
$\euP_2^{\lceil i/2\rceil }/\euP_1^{\lceil i/4\rceil }$.  We can
change $\rho_m$ by an element in $\euP_1^{\lceil i/4\rceil }$ and
still have a $\euO_T$-basis. 
So when $\rho_m-(\sigma +1)\alpha_m\in \euP_1^{\lceil i/4\rceil }$,
the difference between
$\rho_m$ and $(\sigma+1)\alpha_m$ is {\em trivial}.

Since $v_2((\sigma+1)\alpha)=v_2(\rho-(\sigma+1)\alpha)$, checking
$\rho_m-(\sigma +1)\alpha_m\in \euP_1^{\lceil i/4\rceil }$ is
equivalent to checking $v_3((\sigma+1)\alpha)\geq i$.  In Case $A$,
because $b_2+b_1\leq 4e_0$ we find that
$v_3((1/2)\cdot(\sigma+1)(\sigma^2+1)\alpha_m)\leq
v_3((\sigma+1)\alpha_m)$.  Therefore, in $A(3)$ through $A(8)$, we may
replace $\rho_m$ by $(\sigma+1)\alpha_m$.  We refrain from doing so in
$A(8)$ as it may hamper our ability to determine the effect of
$\mbox{Tr}_{3,2}$ on $\overline{\rho}$. We will return to this issue
in \S 3.3.4.  In Case $B$, because $b_2>4e_0-2b_1$ we find
$v_3(\overline{2\alpha})<v_3((\sigma +1)\alpha)$.  We may replace
$\rho$ in $B(3)$ through $B(8)$.  For similar reasons, we refrain in
$B(8)$.  In Cases $C$ and $D$, $b_3> 2b_2+2b_1$ (since $b_3=b_2+4e_0$
and $b_2<4e_0-2b_1$).  As a consequence, $v_3(\overline{(\sigma
+1)(\sigma^2+1)\alpha})<v_3((\sigma +1)\alpha)$. We may replace $\rho$
in $C(3)$ through $C(8)$, and in $D(3)$ through $D(6)$.  In Cases $E$
through $H$, we clearly have $v_3(\alpha)<v_3((\sigma +1)\alpha)$. We
may replace $\rho$ in $E(1)$ or $E(3)$ -- $E(8)$, $F(1)$ or $F(3)$ --
$F(8)$, $G(1)$ or $G(3)$ -- $G(8)$, $H(1)$ or $H(3)$ -- $H(8)$.  We
replace $\rho$ everywhere that we may, {\em except} that we refrain
for
\begin{equation} A(8), B(8),
C(8), D(7), D(8), E(8), F(7),
F(8), G(6), G(7), H(6). 
\end{equation}

Now we consider the difference between $\overline{\rho}$ and $(\sigma +1)\overline{\alpha}$ and replace
$\overline{\rho}$ with $(\sigma +1)\overline{\alpha}$
($\overline{2\rho}$ with $(\sigma +1)\overline{2\alpha}$) in 
\begin{equation}
E(1), F(1), G(1), G(8), H(1), H(7), H(8).
\end{equation}
Since $(\sigma^4+1)\cdot
[\overline{\rho} - (\sigma +1)\overline{\alpha}]=0$, we may use Lem
2.1(2) and find an element $\omega\in K_3$ with $v_3(\omega )=
2b_2+b_1-2b_3$ so that $(\sigma^4-1)\omega = \overline{\rho} - (\sigma
+1)\overline{\alpha}$. As long as $b_3<8e_0-3b_1$, which holds in
Cases $E$ through $H$, we have
$v_3(\overline{\rho})=v_3(\overline{\rho}+2\omega)$.  On the basis of
valuation, we may replace $\overline{\rho}$ with
$\overline{\rho}+2\omega$ and still have a basis ({\em i.e.}
Observation (2)). Now since $(\overline{\rho}+2\omega) - (\sigma
+1)\overline{\alpha}=(\sigma^4+1)\omega\in K_2$, we may replace
$(\overline{\rho}+2\omega)$ with $(\sigma+1)\overline{\alpha}$ and
still have a basis. All we need is
$v_3((\sigma+1)\overline{\alpha})\geq i$. But this clearly holds since
$v_3(\overline{\alpha})\geq i$.

\subsubsection{Nakayama's Lemma and an $\euO_T[G]/\langle\mbox{\rm Tr}_{3,2}\rangle$-basis for $\euP_3^i/\euP_2^{\lceil i/2\rceil }$}
The collection of $\overline{X}$ in our bases provide an
$\euO_T$-basis for $\euP_3^i/\euP_2^{\lceil i/2\rceil}$.  As in \S
2.2.3, whenever $\overline{X}$ and $(1/2)\cdot X$ appear in the same
row, we may replace $\overline{X}$ with $\overline{X}-(1/2)\cdot X$
and still have a $\euO_T$-basis.  Since
$\mbox{Tr}_{3,2}(\overline{X}-(1/2)\cdot X)=0$, we relabel and assume,
without loss of generality, that for these $\overline{X}$'s,
$\mbox{Tr}_{3,2}\overline{X}=0$. Let $\mathcal{T}_{=0}$ denote this
set (trace zero).  Let $\mathcal{T}_{\neq 0}$ denote the
set of $\overline{X}$'s with $X$ in the same row. For each such
$\overline{X}\in \mathcal{T}_{\neq 0}$,
$\mbox{Tr}_{3,2}\overline{X}\not\equiv 0\bmod 2$. This is the set of
trace {\em not} zero.  Note that $\mbox{Tr}_{3,2}\mathcal{T}_{\neq 0}$
is an $\euO_T/2\euO_T$-basis for
$\mbox{Tr}_{3,2}\euP_3^i/2\euP_2^{\lceil i/2\rceil}$.  Following \S
2.2.3, we select from $\mathcal{T}_{\neq 0}$ a set $\mathcal{S}$
(notation as in \S 2.2.3) such that $\mbox{Tr}_{3,2}\mathcal{S}$ is a
$\euO_T/2\euO_T$-basis for $\mbox{Tr}_{3,2}\euP_3^i/ (
(\sigma-1)\mbox{Tr}_{3,2}\euP_3^i+2\euP_2^{\lceil i/2\rceil} )$. It
turns out that just as in \S 2.2.3, $\mathcal{S}$ is the set of
left-most $\overline{X}$ for which $X$ appears in the same row, {\em
except} that $\mathcal{S}$ contains both $\overline{X}$'s in
$\mathcal{T}_{\neq 0}$ from $B(3)$, $C(3)$, $D(3)$.

Note that $\sigma$ acts trivially (modulo $2$) upon $(\sigma
+1)(\sigma^2+1)\alpha$ and $2\alpha$ in $B(3)$, $C(3)$ and $D(3)$.
These elements are linearly independent over
$\euO_T/2\euO_T[G]$. Since both contribute to the
$\euO_T/2\euO_T$-basis for $\mbox{Tr}_{3,2}\euP_3^i/2\euP_2^{\lceil
i/2\rceil}$, both $\overline{(\sigma +1)(\sigma^2+1)\alpha}$ and
$\overline{2\alpha}$ are in $\mathcal{S}$.  When a row contributes
exactly one $\overline{X}$ to $\mathcal{T}_{\neq 0}$, the phrase
`left--most' is unnecessary. Indeed $\sigma$ acts trivially (modulo
$2$) on the lone $X=\mbox{Tr}_{3,2}\overline{X}$, and since $X$ is
needed for the $\euO_T/2\euO_T$-basis for
$\mbox{Tr}_{3,2}\euP_3^i/2\euP_2^{\lceil i/2\rceil}$, $\overline{X}$
must appear in $\mathcal{S}$. Note this is the only situation to
consider in Case $A$.

In the other cases, we need to show that each $X$, corresponding to the
left--most $\overline{X}$ of $\mathcal{T}_{\neq 0}$, generates over
$\euO_T/2\euO_T[G]$ all other elements in the same row (in
$\mbox{Tr}_{3,2}\mathcal{T}_{\neq 0}$). This is easy to see for rows
$E(1)$, $E(5)$, $F(1)$, $F(5)$, $G(1)$, $G(5)$, $G(8)$ and $H(1)$,
$H(5)$, $H(7)$, $H(8)$.  More work is required for rows $D(7)$,
$F(7)$, $G(6)$, $G(7)$, $H(6)$.  Note that
$\rho-(\sigma+1)\alpha_m=(\sigma^2+1)\alpha_{m-t}$ or
$(\sigma+1)(\sigma^2+1)\alpha_{m-s}$ depending upon $b_2>3b_1$ or
$b_2=3b_1$, respectively.  If $\rho-(\sigma+1)\alpha_m=
(\sigma+1)(\sigma^2+1)\alpha_{m-s}$, then $(\sigma
-1)\rho=(\sigma^2+1)\alpha -2\alpha\equiv (\sigma^2+1)\alpha \bmod
2\euP_2^{\lceil i/2\rceil}$.  So $\rho$ generates
$(\sigma^2+1)\alpha$.  If $\rho-(\sigma+1)\alpha_m=
(\sigma^2+1)\alpha_{m-t}$ the analysis is a little more involved.
Note $(\sigma-1)\rho-(\sigma^2+1)\alpha \equiv
(\sigma-1)(\sigma^2+1)\alpha_{m-t}\bmod 2\euP_2^{\lceil i/2\rceil}$.
For $m$ associated with $D(7)$, $F(7)$, $G(6)$, $G(7)$, $H(6)$, check
that $m-t$ lies in $D(3)$, $F(4)$, $G(4)$, $H(4)$ or later. In any
case $\overline{(\sigma+1)(\sigma^2+1)\alpha_{m-t}}\in\euP_3^i$.  So
$\overline{\rho}_m$ and another $\overline{X}$, namely
$\overline{(\sigma +1)(\sigma^2+1)\alpha}_{m-t}$, combine together to
generate $(\sigma^2+1)\alpha_m$.

Apply Lemma 2.2 and
extend $\mathcal{S}$ to an
$\euO_T[G]/\langle\mbox{Tr}_{3,2}\rangle$-basis for
$\euP_3^i/\euP_2^{\lceil i/2\rceil}$.  Except in Cases $B$,
$C$, $D$ (where a row contributes more than one element), we may
assume that this basis is the set of left--most elements
$\overline{X}$, one from each row.

\subsubsection{Essentially Trivial Difference} In \S 3.3.2
we did not replace $\rho$ by $(\sigma +1)\alpha$ in rows $A(1)$,
$A(2)$, $B(1)$, $B(2)$, $C(1)$, $C(2)$, $D(1)$, $D(2)$, $E(2)$,
$F(2)$, $G(2)$, $H(2)$.  It was not clear that the difference
$\rho-(\sigma +1)\alpha$ lay in $\euP_1^{\lceil i/4\rceil}$.  Neither
did we replace $\rho$ by $(\sigma +1)\alpha$ in the rows listed in
(3.3).  In this section we remedy this situation. We
show, except in four cases, $C(2),
D(2), E(2), F(2)$, we may change our basis so that the Galois
action is {\bf\em as if} $\rho$ had been replaced by $(\sigma
+1)\alpha$ ($\overline{\rho}$ by $(\sigma +1)\overline{\alpha}$).

We begin with {\em Case} $A$, explaining why the difference between
$\rho$ and $(\sigma +1)\alpha$ is {\em essentially trivial} and then
determine the Galois module structure (to illustrate the process).
Consider $A(1), A(2)$ and $A(8)$.  Recall  there are three
expressions for $\rho_m$ corresponding to $3b_1<b_2<4e_0-b_1$,
$b_2=3b_1$, and $b_2=4e_0-b_1$. Suppose $3b_1<b_2<4e_0-b_1$, and 
$\rho_m=(\sigma+1)\alpha_m + (\sigma^2-1)\alpha_{m-t}$.  Consider
$\rho_m$ in $A(8)$.  Since $b_1+b_2< 4e_0$,
$v_3(\overline{\rho}_m)\leq v_3(\overline{2\alpha}_{m-t})$. So for $m$
in $A(8)$, $m-t$ is in $A(4)$ or later.  In any case,
$(\sigma-1)\alpha_{m-t}=(\sigma
+1)\alpha_{m-t}-2\alpha_{m-t}\in\euP_3^i$ and
$(1/2)(\sigma-1)(\sigma^2+1)\alpha_{m-t}=
(1/2)(\sigma+1)(\sigma^2+1)\alpha_{m-t}-(\sigma^2+1)\alpha_{m-t}\in\euP_3^i$
({\em i.e.} these elements are available).  We replace $\alpha_m$ with
$x=\alpha_m+(\sigma-1)\alpha_{m-t}
-(1/2)(\sigma-1)(\sigma^2+1)\alpha_{m-t}$.  Note $(\sigma
+1)x=\rho$ and $(\sigma ^2+1)x=(\sigma^2+1)\alpha_m$.  The Galois
action on $x$ and $\rho_m$ is the same as the Galois action on
$\alpha_m$ and $(\sigma+1)\alpha_m$. It is {\bf\em as if} $\rho_m$ had
been replaced by $(\sigma+1)\alpha_m$ and $\overline{\rho}_m$ by
$(\sigma+1)\overline{\alpha}_m$.  Now consider $A(1)$ and $A(2)$,
$\rho_m=(\sigma+1)\alpha_m +(1/2)(\sigma^2-1)\alpha_{m-t+e_0}$.  Since
$v_2(\rho_m)<v_2((1/2)(\sigma +1)(\sigma^2+1)\alpha_{m-t+e_0})$, for
$m$ in $A(1)$ or $A(2)$, $m-t+e_0$ lies in $A(3)$ or later. In any
case, $(1/2)(\sigma -1)(\sigma^2-1)\alpha_{m-t+e_0}$ is available. So
in $A(1)$ and $A(2)$, we replace $2\alpha_m$ by $2\alpha_m
-(1/2)(\sigma -1)(\sigma^2-1)\alpha_{m-t+e_0}$.  The effect of this
replacement on the Galois action is, again, the same {\bf\em as if} we
replaced $\rho_m$ by $(\sigma +1)\alpha_m$.

Now suppose $b_2=3b_1$ and $\rho_m=(\sigma+1)\alpha_m +
(\sigma+1)(\sigma^2+1)\alpha_{m-s}$. Note $s=b_1$.  Starting with the
smallest $m$ such that $i\leq v_3(\overline{\rho}_m)$ we replace
$\alpha_m$ by $\alpha_m+(1/2)(\sigma+1)\alpha_{m+e_0-b_1}$ so long as
$m+e_0-b_1$ is associated with $A(8)$. If $i\leq v_3(\rho_{m-b_1})$,
we replace $\alpha_m$ by $\alpha_m+(\sigma^2+1)\alpha_{m-b_1}$. In any
case, we can systematically replace $\alpha_m$ by
$x=\alpha_m+(1/2)(\sigma^2+1)\alpha_{m+e_0-b_1}$ or
$\alpha_m+(\sigma^2+1)\alpha_{m-b_1}$(1/2)$(\sigma ^2+1)\alpha_m$ by
$(1/2)(\sigma ^2+1)x$ and $(1/2)(\sigma+1)(\sigma ^2+1)\alpha_m$ by
$(1/2)(\sigma +1)(\sigma ^2+1)x$. 
The Galois action after this change of basis is {\bf\em as if}
$\overline{\rho}_m=(\sigma+1)\overline{\alpha}_m$ and $\rho_m=(\sigma
+1)\alpha_m$. Consider $A(1)$ and $A(2)$. Note  $(\sigma
-1)\rho_m=(\sigma -1)\cdot (\sigma +1)\alpha_m$. Moreover, for $m$
associated with these two cases,
$(\sigma+1)(\sigma^2+1)\alpha_{m+e_0-b_1}$ and
$(\sigma^2+1)\alpha_{m+e_0-b_1}$ are available elsewhere in our
basis. So we replace $(\sigma^2+1)\alpha_m$ by
$(\sigma^2+1)(\alpha_m+\alpha_{m+e_0-b_1})$ and $(\sigma
+1)(\sigma^2+1)\alpha_m$ by $(\sigma
+1)(\sigma^2+1)(\alpha_m+\alpha_{m+e_0-b_1})$. Note 
for $m$ associated with $A(2)$, $m+e_0-b_1$ is associated with $A(3)$
or later. We achieve the desired effect by replacing
$\overline{(\sigma
+1)(\sigma^2+1)\alpha}_m$ with $\overline{(\sigma
+1)(\sigma^2+1)\alpha}_m+\overline{(\sigma
+1)(\sigma^2+1)\alpha}_{m+e_0-b_1}$.

This leaves $b_2=4e_0-b_1$. Because this case is more complicated
(recall Remark 3.5: $\rho_m$ is `torn' between $\alpha_m$ and
$\alpha_{m+e_0-b_1}$), we first determine the Galois module structure
for $b_2<4e_0-b_1$.  Each $m$ in $A(1)$ results in an
$\euO_T\otimes_{\mathbb{Z}_2}(\mathcal{R}_3\oplus\mathcal{H})$; $m$ in $A(2)$
in an $\euO_T\otimes_{\mathbb{Z}_2}\mathcal{H}_2$; $m$ in $A(3)$ in an
$\euO_T\otimes_{\mathbb{Z}_2}(\mathcal{R}_3\oplus\mathcal{M})$; $m$ in $A(4)$
in an $\euO_T\otimes_{\mathbb{Z}_2}\mathcal{M}_1$; $m$ in $A(5)$ in an
$\euO_T\otimes_{\mathbb{Z}_2}(\mathcal{R}_3\oplus\mathcal{L})$; $m$ in $A(6)$
in an $\euO_T\otimes_{\mathbb{Z}_2}\mathcal{L}_3$; $m$ in $A(7)$ in an
$\euO_T\otimes_{\mathbb{Z}_2}(\mathcal{R}_3\oplus\mathcal{I})$; $m$ in $A(8)$
in an $\euO_T\otimes_{\mathbb{Z}_2}\mathcal{I}_2$.  Counting the number of
$m$ associated with each $A(j)$ yields the first column of Table 2.

Now consider $b_2=4e_0-b_1$. Because $v_2(\rho_m)=2b_2-b_1+4m$,
the number of $m$ associated with $A(1)$ and $A(7)$ are
different.  The number for $A(7)$ is $e_0-b_1$ too low, while 
$A(1)$ is $e_0-b_1$ too high.  We seem to be missing
$e_0-b_1$ of $\euO_T\otimes_{\mathbb{Z}_2}\mathcal{I}$ and have $e_0-b_1$ too
many of $\euO_T\otimes_{\mathbb{Z}_2}\mathcal{H}$.  Let us look at this more
carefully.  Note  $\overline{\rho}_m$ in $A(8)$ maps (via
$\mbox{Tr}_{3,2}$) to
$$\rho_m=(\sigma +1)(\alpha_m-(1/2)(\sigma^2+1)\alpha_m) + 
\begin{cases}
(1/2)(\sigma+1)(\sigma^2+1)\alpha_{m+e_0-b_1} &  \\
(\sigma+1)(\sigma^2+1)\alpha_{m-b_1} & 
\end{cases}
$$
So $\overline{\rho}_m$ maps into the $\euO_T$-module spanned by
$\alpha_m-(1/2)(\sigma^2+1)\alpha_m$ and $(\sigma
+1)(\alpha_m-(1/2)(\sigma^2+1)\alpha_m)$ along with {\em either}
$(1/2)(\sigma^2+1)\alpha_{m+e_0-b_1}$ and
$(1/2)(\sigma+1)(\sigma^2+1)\alpha_{m+e_0-b_1}$ {\em or} 
$(\sigma^2+1)\alpha_{m-b_1}$ and
$(\sigma+1)(\sigma^2+1)\alpha_{m-b_1}$. In any case, the elements
$(1/2)(\sigma^2+1)\alpha_m$ and $(1/2)(\sigma+1)(\sigma^2+1)\alpha_m$
for $\lceil (i+b_3-4b_2+2b_1)/8\rceil \leq m \leq \lceil
(i+b_3-4b_2+2b_1)/8\rceil + e_0-b_1-1$ are not associated with a
$\overline{\rho}_m$ in $A(8)$.  The $\rho_m$ in $A(1)$ map
to $(\sigma^2-1)\alpha_m$ (under $(\sigma -1)$) and so
$(\sigma+1)(\sigma^2+1)\alpha_{m+e_0-b_1}$ (under $(\sigma^2+1)$)
yielding a $\mathcal{H}$, unless $m+e_0-b_1$ is associated with
$A(2)$. In fact, there are $e_0-b_1$ $\rho_m$ that map into $A(2)$
under $(\sigma^2+1)$. For each $m$ in $A(2)$ we have
$(\sigma^4+1)\overline{(\sigma+1)(\sigma^2+1)\alpha}_m
=(\sigma^2+1)\rho_{m-e_0+b_1} =(\sigma+1)(\sigma^2+1)\alpha_m$,
yielding a copy of $\mathcal{H}_2$. But for the last $e_0-b_1$
elements $\rho_m$ in $A(2)$, namely those $m$ such that $m+e_0-b_1$ is
in $A(3)$ we may replace $\rho_m$ by
$\rho_m-(1/2)(\sigma+1)(\sigma^2+1)\alpha_{m+e_0-b_1}$.  For each
of these $m$ we have the $\euO_T[G]$-submodule spanned by
$\rho_m-(1/2)(\sigma+1)(\sigma^2+1)\alpha_{m+e_0-b_1}$ and
$(\sigma^2-1)\alpha_m$.  These $e_0-b_1$ together with the elements
left out of a module in $A(8)$ yield a $e_0-b_1$ copies of
$\mathcal{I}$, precisely making up the counts.

\vspace*{1mm}

\noindent{\em Cases} $B$ -- $H$: In the remaining cases, we only have
two situations: $b_2=3b_1$ and $3b_1<b_2<4e_0-b_1$. Consider
$3b_1<b_2<4e_0-b_1$ first, and $\rho_m=(\sigma+1)\alpha_m
+(\sigma^2\pm 1)\alpha_{m-t}$ where we may choose between $\pm$ as we
like.  We are concerned with the image of the trace,
$\mbox{Tr}_{3,2}$, in particular
$\mbox{Tr}_{3,2}\overline{\rho}_m=(\sigma+1)\alpha_m+(\sigma^2 +
1)\alpha_{m-t}$, for $\overline{\rho}_m$ appearing in $B(8)$, $C(8)$,
$D(7)$, $D(8)$, $E(8)$, $F(7)$, $F(8)$, $G(6)$, $G(7)$, and $H(6)$.
Note  if $\overline{(\sigma^2+1)\alpha_{m-t}}\in\euP_3^i$, we
may replace $\overline{\rho}_m$ with $\overline{\rho}_m -
\overline{(\sigma^2+1)\alpha_{m-t}}$. So if 
$(\sigma^2+1)\alpha_{m-t}$ appears in $B(6)$, $C(6)$, $D(6)$, $E(5)$,
$F(5)$, $G(5)$, $H(5)$ or later we may replace $\rho_m$ with
$(\sigma+1)\alpha_m$ and $\overline{\rho}_m$ with $\overline{\rho}_m -
\overline{(\sigma^2+1)\alpha_{m-t}}$. The later replacement exhibits
the same Galois action as a replacement of $\overline{\rho_m}$ by
$(\sigma+1)\overline{\alpha}_m$. Without loss of generality we will
call it a replacement of $\overline{\rho_m}$ by
$(\sigma+1)\overline{\alpha}_m$.

Since $b_2\leq 4e_0-b_1$, $v_3(\overline{2\alpha}_m)\geq
v_3(\overline{\rho}_m)$.  What happens when $(\sigma^2+1)\alpha_{m-t}$
appears in $B(3)$ -- $B(5)$, $C(3)$ -- $C(5)$, $D(3)$ -- $D(5)$,
$E(4)$, $F(6)$, $G(4)$, $H(6)$?  In this case $(\sigma
-1)\alpha_{m-t}= (\sigma +1)\alpha_{m-t}-2\alpha_{m-t}\in \euP_3^i$.
In $B(8)$, $C(8)$, $D(7)$, $D(8)$, $E(8)$, $F(8)$, $G(6)$, $G(7)$,
$H(6)$, we replace $\alpha_m$ with $\alpha_m+ (\sigma
-1)\alpha_{m-t}$, and $(\sigma^2+1)\alpha_m$ with
$(\sigma^2+1)\alpha_m+ (\sigma -1)(\sigma^2+1)\alpha_{m-t}$.  Note
$\rho_m=(\sigma+1)\cdot [ \alpha_m+(\sigma-1)\alpha_{m-t}]$.  The
Galois action upon these basis elements:
$\mbox{Tr}_{3,2}\overline{\rho}_m=\rho_m=(\sigma+1)\cdot [
\alpha_m+(\sigma-1)\alpha_{m-t}]$, $(\sigma^2+1)\cdot [
\alpha_m+(\sigma-1)\alpha_{m-t}]= (\sigma^2+1)\alpha_m+ (\sigma
-1)(\sigma^2+1)\alpha_{m-t}$, and
$(\sigma+1)\cdot[(\sigma^2+1)\alpha_m+ (\sigma
-1)(\sigma^2+1)\alpha_{m-t}] =(\sigma^2+1)\rho_m= (\sigma
+1)(\sigma^2+1)\alpha_m$, is similar to the Galois action upon:
$(\sigma+1)\overline{\alpha}_m$, $\alpha_m$, $(\sigma+1)\alpha_m$,
$(\sigma^2+1)\alpha_m$, $(\sigma+1)(\sigma^2+1)\alpha_m$.  We may
assume $(\sigma+1)\overline{\alpha}_m$ and $(\sigma+1)\alpha_m$ appear
instead of $\overline{\rho}_m$ and $\rho_m$.

Now consider the appearance of $\rho$ in $B(1)$,
$B(2)$, $C(1)$, $D(1)$, $G(2)$, $H(2)$.  Suppose 
$\rho_m=(\sigma+1)\alpha_m+(1/2)\cdot(\sigma^2-1)\alpha_{m+e_0-t}$.
One may check  $v_3(\rho_m)\leq
v_3((1/2)(\sigma+1)(\sigma^2+1)\alpha_{m+e_0-t}$ and 
$v_3(\overline{2\rho}_{m+e_0-t})\leq v_3(\overline{4\alpha_m})$.  So
$(1/2)(\sigma+1)(\sigma^2+1)\alpha_{m+e_0-t}$ appears in $B(4)$ --
$B(7)$, $C(6)$ -- $C(8)$ or $D(6)$ -- $D(8)$. Note  in these sets
of elements, $\rho_{m+e_0-t}$ has already been replaced by
$(\sigma+1)\alpha_{m+e_0-t}$. Importantly,
$(1/2)(\sigma-1)(\sigma^2+1)\alpha_{m+e_0-t}$ along with
$(\sigma-1)\alpha_{m+e_0-t}$ are available to us. We replace
$2\alpha_m$ with
$2\alpha_m-(1/2)(\sigma-1)(\sigma^2+1)\alpha_{m+e_0-t}+
(\sigma-1)\alpha_{m+e_0-t}=2\alpha_m-(1/2)(\sigma-1)(\sigma^2-1)\alpha_{m+e_0-t}$
in $B(1)$, $B(2)$, $C(1)$ and $D(1)$. The effect of this change of
basis is the same as if we replaced $\rho_m$ by $(\sigma+1)\alpha_m$.

Now consider $G(2)$ and $H(2)$. Again
$\rho_m=(\sigma+1)\alpha_m+(1/2)(\sigma^2-1)\alpha_{m+e_0-t}$.  In $G$
and $H$, $b_3\leq 8e_0-2b_2$. As a result, $v_3(\rho_m)\leq
v_3(\overline{(\sigma -1)\alpha_{m+e_0-t}})$.  Note  we refer to
$\overline{(\sigma -1)\alpha_{m+e_0-t}}$ and not $(\sigma
-1)\overline{\alpha_{m+e_0-t}}$. The valuation of the first is $b_1$
more than the valuation of the second.  As one may check
$v_3(\rho_m)\leq v_3((1/2)(\sigma+1)(\sigma^2+1)\alpha_{m+e_0-t})$, so
$(1/2)(\sigma+1)(\sigma^2+1)\alpha_{m+e_0-t}$ appears in $G(7)$,
$G(8)$ or $H(8)$.  If $(\sigma+1)(\sigma^2+1)\alpha$ appeared in
$G(1)$ or $H(1)$, $(\sigma^2+1)\alpha_{m-t}$ would be available and so
we would replace $\rho_m$ with $\rho_m-(\sigma^2+1)\alpha_{m-t}$.  If
$(1/2)(\sigma+1)(\sigma^2+1)\alpha_{m+e_0-t}$ appears in $G(7)$, then
we may assume  $\overline{(\sigma -1)\alpha_{m+e_0-t}}$ appears
there instead of $\overline{\rho}_{m+e_0-t}$, because
$v_3(\overline{(\sigma^2 +1)\alpha_{m+e_0-2t}})=v_3(\overline{(\sigma
-1)\alpha_{m+e_0-2t}})\geq i$, and we would have replaced
$\overline{\rho}_{m+e_0-t}$ previously in our discussion with
$\overline{\rho}_{m+e_0-t}-\overline{(\sigma^2 +1)\alpha_{m+e_0-2t}}$.
We may now replace $\overline{2\alpha}_m$ with $\overline{2\alpha}_m-
\overline{(\sigma -1)\alpha_{m+e_0-t}}$. We replace
$(\sigma^2+1)\alpha_m$ with $(\sigma^2+1)\alpha_m + (1/2)(\sigma
-1)(\sigma^2)\alpha_{m+e_0-t}$.  We may assume without loss of
generality that $(\sigma+1)\alpha_m$ appears in $G(2)$ and $H(2)$
instead of $\rho_m$.


Now we work with Cases $B$ through $H$ under the assumption 
$b_2=3b_1$.  So $\rho_m=(\sigma+1)\cdot [
\alpha_m+(\sigma^2+1)\alpha_{m-b_1}]$.  First note  if
$(\sigma+1)(\sigma^2+1)\alpha_{m-b_1}$ appears in $B(2)$, $C(3)$,
$D(3)$, $E(4)$, $F(4)$, $G(4)$, $H(4)$, or later we may replace
$\overline{\rho}_m$ in $B(8)$, $C(8)$, $D(7)$, $D(8)$, $E(7)$, $F(7)$,
$G(5)$, $G(6)$, $H(5)$ with
$\overline{\rho}_m-\overline{(\sigma+1)(\sigma^2+1)\alpha}_m$.
Suppose $(\sigma^2+1)\alpha_{m-b_1}$ appears elsewhere.  In $B$, these
elements can appear in $B(1)$, $B(2)$, or as
$(1/2)\cdot(\sigma^2+1)\alpha_{m+e_0-b_1}$ elsewhere in $B(8)$.  In
cases $C$ through $H$, since $b_1<4e_0/5$,
$v_3(\overline{\rho}_m)\leq v_3(\rho_{m-b_1})$. So
$(\sigma^2+1)\alpha_{m-b_1}$ appears in $C(1)$, $C(2)$, $D(1)$,
$D(2)$, $E(2)$, $E(3)$, $F(2)$, $F(3)$, $G(2)$, $G(3)$, $H(2)$,
$H(3)$. In these cases, we may either replace $\alpha_m$ with
$\alpha_m+(1/2)\cdot(\sigma^2+1)\alpha_{m+e_0-b_1}$ or
$\alpha_m+(\sigma^2+1)\alpha_{m+-b_1}$.
If for example, we replace $\alpha_m$ with
$\alpha_m+(\sigma^2+1)\alpha_{m-b_1}$, $(\sigma^2+1)\alpha_m$ with
$(\sigma^2+1)\alpha_m+2(\sigma^2+1)\alpha_{m-b_1}$, and
$(\sigma+1)(\sigma^2+1)\alpha_m$ with
$(\sigma+1)(\sigma^2+1)\alpha_m+2(\sigma+1)(\sigma^2+1)\alpha_{m-b_1}$,
then the Galois action on this new basis is the same {\em as if}
$(\sigma+1)\overline{\alpha}_m$ and $(\sigma+1)\alpha_m$ appear
instead of $\overline{\rho}_m$ and $\rho_m$.

We now concern ourselves with $B(1)$, $B(2)$, $C(1)$, $D(1)$, $G(2)$
and $H(2)$.  Check  $v_3((\sigma^2+1)\alpha_{m+e_0-b_1})\geq
v_3(\rho_m)$. We replace $(\sigma^2+1)\alpha_m$ with
$(\sigma^2+1)\alpha_m + (\sigma^2+1)\alpha_{m+e_0-b_1}$, and
$(\sigma+1)(\sigma^2+1)\alpha_m$ with $(\sigma+1)(\sigma^2+1)\alpha_m
+ (\sigma+1)(\sigma^2+1)\alpha_{m+e_0-b_1}$.  In $B(2)$,
$v_3(\overline{(\sigma+1)(\sigma^2+1)\alpha}_{m+e_0-b_1})\geq
v_3(\overline{(\sigma+1)(\sigma^2+1)\alpha}_m)$, we replace
$(\sigma+1)(\sigma^2+1)\alpha_m$ with $(\sigma+1)(\sigma^2+1)\alpha_m+
(\sigma+1)(\sigma^2+1)\alpha_{m+e_0-b_1}$.  All this has the same
effect upon the Galois action as a replacement of $\rho_m$ by
$(\sigma+1)\alpha_m$. 

\subsubsection{Non-Trivial Difference}
We consider $\rho$ in $C(2)$,
$D(2)$, $E(2)$, $F(2)$. 

First consider the case $b_2= 3b_1$ where
$\rho_m=(\sigma+1)\alpha_m+
(1/2)(\sigma+1)(\sigma^2+1)\alpha_{m+e_0-b_1}$.  Note  $C$ and $E$
do not intersect the line $b_2=3b_1$. We focus on $D(2)$, $F(2)$.
In $D$ with $b_2=3b_1$, we have $b_1<4e_0/5$. So
$v_3(\overline{2\alpha}_m)\leq v_3(\alpha_{m+e_0-b_1})$. Since
$v_3(\overline{2(\sigma+1)(\sigma^2+1)\alpha}_m) \leq
v_3((\sigma+1)(\sigma^2+1)\alpha_{m+e_0-b_1})$, for $m$ associated
with $D(2)$, $(\sigma+1)(\sigma^2+1)\alpha_{m+e_0-b_1}$ appears in
$D(4)$, or $(1/2)(\sigma+1)(\sigma^2+1)\alpha_{m+e_0-b_1}$ appears in
$D(5)$ or later. If $(1/2)(\sigma+1)(\sigma^2+1)\alpha_{m+e_0-b_1}$ is
available, we may replace $\rho_m$ with $(\sigma+1)\alpha_m$.  The
Galois action when $m$ is in $D(2)$ and $m+e_0-b_1$ is in $D(4)$ is
our primary concern. But first consider $F$ (or $\overline{F}$) with
$b_2=3b_1$. Note then  $b_3\leq 8e_0+2b_2-8b_1$. So
$v_3(\rho_m)\leq
v_3(\overline{(\sigma+1)(\sigma^2+1)\alpha}_{m+e_0-b_1})$.  Since
$b_3\leq 8e_0+2b_2-8b_1$, $v_3(\alpha_m)\leq
v_3(\overline{2(\sigma^2+1)\alpha}_{m+e_0-b_1})$.  So for $m$
associated with $F(2)$, $(\sigma+1)(\sigma^2+1)\alpha_{m+e_0-b_1}$
appears in $F(4)$, or in $F(5)$ or later. If $m+e_0-b_1$ is associated
with $F(5)$ or later, we have
$\overline{(\sigma^2+1)\alpha}_{m+e_0-b_1}$ available.  We replace
$\overline{2\alpha}_m$ with $\overline{2\alpha}_m+
\overline{(\sigma^2+1)\alpha}_{m+e_0-b_1}$.  We replace $(\sigma^2
+1)\alpha_m$ and $(\sigma +1)(\sigma^2 +1)\alpha_m$ with $(\sigma^2
+1)\alpha_m+(\sigma^2 +1)\alpha_{m+e_0-b_1} $ and $(\sigma
+1)(\sigma^2 +1)\alpha_m+(\sigma +1)(\sigma^2 +1)\alpha_{m+e_0-b_1}$.
The effect of these changes upon the Galois action is the same as the
replacement of $\rho_m$ by $(\sigma+1)\alpha_m$.  This leaves the
situation when $m$ belongs to $D(2)$, $F(2)$ while $m+e_0-b_1$ belongs
to $D(4)$, $F(4)$.  In both of these cases, we replace
$\overline{(\sigma+1)(\sigma^2+1)\alpha}_{m+e_0-b_1}$ with
$\overline{(\sigma+1)(\sigma^2+1)\alpha}_{m+e_0-b_1} +
(\sigma+1)\overline{2\alpha}_m-\rho_m$. This new basis element has
trace, $\mbox{Tr}_{3,2}$, zero. For each such pair $(m,m+e_0-t)$ we
get a copy of $\mathcal{H}_1\mathcal{G}\oplus \mathcal{R}_3$.

Let us now turn to the case where $3b_1<b_2<4e_0-b_1$ and 
$\rho_m=(\sigma+1)\alpha_m+(1/2)(\sigma^2+1)\alpha_{m+e_0-t}$.
Consider cases $C$ and $E$. Because
$v_3(\overline{2\alpha})\leq v_3(\overline{(\sigma ^2+1)\alpha})$, if
$m$ appears in $C(2)$, then $m+e_0-t$ appears in $C(6)$ or later. Since
$v_3((\sigma^2+1)\alpha_{m+e_0-t}>
v_3(\overline{2(\sigma+1)(\sigma ^2+1)\alpha})$, not every $m+e_0-t$
is in $C(6)$ when $m$ is in $C(2)$. Since $v_3(\rho)\leq
v_3((1/2)(\sigma ^2+1)\alpha)$, if $m$ appears in $E(2)$, then
$m+e_0-t$ appears in $E(6)$ or later. Since
$v_3((\sigma^2+1)\alpha_{m+e_0-t}> v_3(2\alpha )$, some
$m+e_0-t$ spill over into $C(7)$.
Consequently, whenever a pair $(m,m+e_0-t)$ has $m$ in $C(2)$, $E(2)$ while
$m+e_0-t$ is in $C(6)$, $E(6)$ we get a copy of 
$\mathcal{H}_1\mathcal{L}\oplus \mathcal{R}_3$. 

Consider cases $D$ and $F$ (including $\overline{F}$). Consider $D$
first.  Since
$v_3(\overline{2\alpha}_m)<v_3(\overline{(\sigma^2+1)\alpha}_{m+e_0-t})$,
for $m$ in $D(2)$, $m+e_0-t$ lands in $D(6)$ or later. Note 
since $v_3(\overline{2\alpha}_m)>v_3(\rho_{m+e_0-t})$, some $m+e_0-t$
land in $D(6)$. Since
$v_3((\sigma^2+1)\alpha_{m+e_0-t})>v_3(\overline{2(\sigma^2+1)\alpha}_m)$,
the collection of $m+e_0-t$ overlap into $D(8)$. When $m+e_0-t$ is in
$D(8)$, the element $(1/2)(\sigma^2+1)\alpha_{m+e_0-t}$ is available
and we replace $\rho_m$ by
$\rho_m-(1/2)(\sigma^2+1)\alpha_{m+e_0-t}=(\sigma +1)\alpha_m$.  For
each pair $(m,m+e_0-t)$ such that $m$ is associated with $D(2)$ and
$m+e_0-t$ is associated with $D(6)$, we get a copy of
$\mathcal{H}_1\mathcal{L}\oplus \mathcal{R}_3$.  What we are
principally concerned with is what happens when for $m$ in
$D(2)$, $m+e_0-t$ is in $D(7)$.  In this case, because
$\rho_{m+e_0-t}=(\sigma+1)\alpha_{m+e_0-t} +
(\sigma^2+1)\alpha_{m+e_0-2t}$, there is some new
interaction to consider.

Suppose $m$ is in $D(2)$, while $m+e_0-t$ is in $D(7)$.  Since
$v_3(\overline{\rho}_{m+e_0-t})\leq v_3(\alpha_{m+e_0-2t})$ and
$v_3(\overline{2(\sigma+1)(\sigma^2+1)\alpha}_m)\leq
v_3((\sigma+1)(\sigma^2+1)\alpha_{m+e_0-2t})$, for $m$ in $D(2)$ and
$m+e_0-t$ in $D(7)$, we find $m+e_0-2t$ is associated with $D(4)$, or
$D(5)$ or later.  Consider $m$ in $D(2)$, $m+e_0-t$ in $D(7)$, and
$m+e_0-2t$ in $D(4)$.  Perform change of basis: Replace
$\overline{2\alpha}_m$ with $\overline{2\alpha}_m
+\overline{2\alpha}_{m+e_0-t} - \overline{2\alpha}_{m+e_0-t}$,
$\rho_m$ with $\rho_m-\alpha_{m+e_0-t}$, $(\sigma^2+1)\alpha_m$ with
$(\sigma^2+1)\alpha_m + (\sigma^2+1)\alpha_{m+e_0-2t}+
1/2(\sigma-1)(\sigma^2+1)\alpha_{m+e_0-t}$, and
$(\sigma+1)(\sigma^2+1)\alpha_m$ with $(\sigma+1)(\sigma^2+1)\alpha_m
+ (\sigma+1)(\sigma^2+1)\alpha_{m+e_0-2t}$.  The effect of these base
changes upon the Galois action is the same as if we were to replace
$\rho_m$ with $(\sigma+1)\alpha_m-(1/2)(\sigma+1)(\sigma^2+1)\alpha_
{m+e_0-2t}$.  Notice the similarity between this expression and the
expression for $\rho_m$ used when $b_2=3b_1$. Consequently, this
scenario results in copies of $\mathcal{H}_1\mathcal{G}\oplus
\mathcal{R}_3$.  (Note if $b_2=3b_1$, then $2t=b_1$.)

In the alternative situation, when $m$ is in $D(2)$, $m+e_0-t$ in $D(7)$,
and $m+e_0-2t$ is in $D(5)$ or later, we perform the same basis
changes. Except, since the element
$(1/2)(\sigma+1)(\sigma^2+1)\alpha_{m+e_0-2t}$ is available, we
replace $\rho_m$ with $\rho_m-\alpha_{m+e_0-t}+
(1/2)(\sigma+1)(\sigma^2+1)\alpha_{m+e_0-2t}$. The effect of this
alternative basis change upon the Galois action is the same as a
simple replacement of $\rho_m$ with $(\sigma+1)\alpha_m$.  We now turn
our attention to Cases $F$ and $\overline{F}$.
Since $0<2b_1$, $v_3(\rho_m)< v_3((1/2)(\sigma +1)(\sigma
^2+1)\alpha_{m+e_0-t}$.  So for $m$ associated with $F(2)$, $m+e_0-t$
is associated with $F(6)$ or later.  We leave it to the reader to
check that $m+e_0-t$ lands in $F(6)$ or $F(7)$.  If $m+e_0-t$ is
associated with $F(7)$, then $m+e_0-2t$ lands in $F(4)$ or $F(5)$.  In
any case, all this is analogous to $D$.

\subsection{The Galois module structure under {\em stable} ramification}
For $p=2$, {\em stable ramification} $b_1\geq e_0$ is nearly {\em
strong ramification} $b_1>(1/2)\cdot pe_0/(p-1)$, (the conditions
differ only when $e_0$ is odd -- $K_0$ tame over $\mathbb{Q}_2$).  In
\cite{elder:bord}, the structure of the ring of integers was
determined under {\em strong ramification} for any prime $p$.  We
revisit that argument extending it to ambiguous ideals and the case
$b_1=e_0$.

Following \S 2.1, $\euP_2^{\lceil i/2\rceil}/\euP_1^{\lceil
i/4\rceil}\cong\left (\euO_T[\sigma ]/\langle \sigma^2+1\rangle\right
)^{e_0}$. So $e_0$ elements generate $\euP_2^{\lceil
i/2\rceil}/\euP_1^{\lceil i/4\rceil}$ over $\euO_T[G]$.  Use Lemmas
3.6, 3.7 to select elements, $\alpha$, with odd valuation $a$ such
that $\lceil i/2\rceil\leq a \leq \lceil i/2\rceil+ 2e_0 -1$.  Each of
these $e_0$ elements gives rise (via the action of $(\sigma \pm 1)$)
to another element, $\rho$ in $K_2$, with odd valuation,
$a+(b_2-b_1)=a+2e_0$.  These $\alpha$ along with their Galois
translates, $\rho\equiv (\sigma \pm 1)\alpha\bmod\euP_1^{\lceil
i/4\rceil} $, have valuations in one--to--one correspondence (via
$v_2$) with the odd integers in $\lceil i/2\rceil,\ldots ,4e_0+\lceil
i/2\rceil-1$, and as a result serve as a $\euO_T$-basis for
$\euP_2^{\lceil i/2\rceil}/\euP_1^{\lceil i/4\rceil}$.  The $\alpha$
provide a $\euO_T[G]/\langle\mbox{Tr}_{2,1}\rangle$-basis.

We need this basis for $\euP_2^{\lceil i/2\rceil}/\euP_1^{\lceil
i/4\rceil}$ to be compatible with our $\euO_T$-basis for
$\euP_1^{\lceil i/4\rceil}$ (as determined as in \S 2.2.1), as well as
our $\euO_T[G]/\langle\mbox{Tr}_{3,2}\rangle$-basis for
$\euP_3^i/\euP_2^{\lceil i/2\rceil}$. First we consider compatibility
with $\euP_1^{\lceil i/4\rceil}$.  The $\euO_T$-basis for
$\euP_1^{\lceil i/4\rceil}$ consists of pairs: either $((\sigma
+1)\eta, \eta)$ or $((\sigma +1)\eta, 2\eta)\in K_0\times K_1$ where
$v_1(\eta)$ is odd.  Because of Lemma 2.1 each coordinate uniquely
determines the other. Now consider pairs where the valuation $v_3$ of
both elements is bound between $i$ and $8e_0+i-1$.  For example, pairs
of the form $((\sigma +1)\eta, \eta)$ appear for $\lceil i/4 \rceil
\leq v_1(\eta) \leq 2e_0+\lceil i/4\rceil -b_1-1$, while pairs of the
form $((\sigma +1)\eta, 2\eta)$ appear for $\lceil i/4 \rceil -b_1
\leq v_1(\eta) \leq \lceil i/4\rceil-1$.  The coordinates of all pairs
provides us with an $\euO_T$ basis for $\euP_1^{\lceil i/4\rceil}$.
Each $\alpha$ with $v_2((\sigma ^2+1)\alpha)\leq 4e_0+\lceil i/2
\rceil -1$ determines (via $(\sigma ^2+1)\alpha\in K_1$) a pair of
elements in the $\euO_T$-basis for $\euP_1^{\lceil i/4\rceil}$.  If
$v_1((\sigma ^2+1)\alpha)$ is odd, then $\alpha$ determines a pair of
the form $((\sigma +1)\eta, 2\eta)$. If even, it determines a pair of
the form $((\sigma +1)\eta, \eta)$.  In general for $\alpha$ with
$v_2((\sigma ^2+1)\alpha)\geq 4e_0+\lceil i/2 \rceil $,
$v_2(1/2(\sigma^2+1)\alpha)\geq \lceil i/2\rceil$. So
$1/2(\sigma^2+1)\alpha$ is available and we may replace $\alpha$ in by
$\alpha-1/2(\sigma^2+1)\alpha$ and still have a basis.  Note $(\sigma
^2+1)\left (\alpha -1/2(\sigma^2+1)\alpha\right )=0$. So we can
assume, without loss of generality, $(\sigma ^2+1)\alpha =0$. This
posses no complication, unless $(\sigma \pm 1)\alpha =\mu +\rho$ with
$\rho$ in the image of $\mbox{Tr}_{3,2}\euP_3^i$. In other words,
$v_2(\rho)\geq \lfloor (b_3+i+1)/2\rfloor$. (Note for $\alpha$ with
$v_2((\sigma ^2+1)\alpha)\leq 4e_0+\lceil i/2 \rceil -1$ and $(\sigma
\pm 1)\alpha =\mu +\rho$, we have $v_2(\rho)< \lfloor
(b_3+i+1)/2\rfloor$.)  For these $\alpha$ (actually $\alpha
-1/2(\sigma^2+1)\alpha$), $\mu$ (actually $\mu -(\sigma \pm
1)1/2(\sigma^2+1)\alpha$) will determine a pair $((\sigma +1)\eta,
2\eta)$ or $((\sigma +1)\eta, \eta)$ in our $\euO_T$-basis for
$\euP_1^{\lceil i/4\rceil}$.  We need simply to show $\mu$ and $\mu
-(\sigma \pm 1)1/2(\sigma^2+1)\alpha$ have the same properties.  We
leave it to the reader to do this (use Lemma 3.6 and 3.7 to show that
the valuations are the same, that $\mu -(\sigma \pm
1)1/2(\sigma^2+1)\alpha\in K_0$ if and only if $\mu\in K_0$).  The
only issue that remains is whether there can be any conflict between a
pair of basis elements for $\euP_1^{\lceil i/4\rceil}$ determined
directly, via $(\sigma ^2+1)\alpha$, and a pair determined indirectly
via $\mu=(\sigma \pm 1)\alpha -\rho$. Note any element in the image of
the trace, $\mbox{Tr}_{2,1}$, has valuation that is larger than the
valuation of every $\mu\in K_1$ that arises from the expression for a
Galois translate $\rho=(\sigma \pm 1)\alpha -\mu$.

We select our $\euO_T[G]$-basis for $\euP_3^i/\euP_2^{\lceil
i/2\rceil}$ now. There is one element $X$ in our
$\euO_T$-basis for $\euP_2^{\lceil i/2\rceil}$ for each valuation
$v_2$ in
\begin{equation}
\lfloor (i+b_3+1)/2\rfloor , \ldots , 4e_0+\lceil i/2 \rceil -1.
\end{equation}
The reader may check  for $v_2(X)$ even, $X=(\sigma^2+1)\alpha$
for some $\alpha$ in our $\euO_T[G]$-basis for $\euP_2^{\lceil
i/2\rceil}/\euP_1^{\lceil i/4\rceil}$.  For $v_2(X)$ odd, since
$\lceil i/2 \rceil + (b_2-b_1) < \lfloor (i+b_3+1)/2\rfloor$,
$X=\rho=(\sigma\pm 1)\alpha - \mu$ also for some $\alpha$.  Use Lem
2.1 to create elements $\overline{X}\in\euP_3^i$ such that
$\mbox{Tr}_{3,2}\overline{X}=X$ and $v_3(\overline{X})=v_3(X)-b_3$.
Note  the elements $(\sigma^2+1)\alpha$ and $\mu$ (from each case)
have expressions in terms of our $\euO_T$-basis
for $\euP_1^{\lceil i/4\rceil}$. These expressions depend solely upon the
valuations of $(\sigma^2+1)\alpha$ and $\mu$.

Before we move on to our result, we should say something about our
basis for $\euP_3^i/\euP_2^{\lceil i/2\rceil}$.
Since $\euO_T[\sigma ]/\langle \sigma^4+1\rangle$ is a principal ideal
domain, $\euP_3^i/\euP_2^{\lceil i/2\rceil}$ is free over
$\euO_T[\sigma ]/\langle \sigma^4+1\rangle$ of rank $e_0$.  Given
elements of $K_2$ with valuation $v_2$ listed in (3.5) we may use Lem
2.1(2) to find elements, $\rho\in \euP_3^i$, whose images under the
trace, $\mbox{Tr}_{3,2}$, lie one--to--one correspondence (via
valuation) with (3.5).  Refer to this set of elements in $\euP_3^i$ as
$\mathcal{S}$.  One can check $b_1+\lfloor (i+b_3+1)/2\rfloor >
4e_0+\lceil i/2 \rceil$. Therefore $(\sigma
-1)\mbox{Tr}_{3,2}\euP_3^i\subseteq 2\euP_2^{\lceil i/2\rceil}$.
Since $\mbox{Tr}_{3,2}\mathcal{S}$ is an $\euO_T$-basis for
$\mbox{Tr}_{3,2}\euP_3^i\subseteq 2\euP_2^{\lceil i/2\rceil}$ and
$\sigma$ acts trivially upon $\mbox{Tr}_{3,2}\euP_3^i\subseteq
2\euP_2^{\lceil i/2\rceil}$ we may use Lemma 2.2 and extend
$\mathcal{S}$ to an $\euO_T[G]/\langle \sigma^4+1\rangle$-basis for
$\euP_3^i/\euP_2^{\lceil i/2\rceil}$.

At this point we may put the preceding discussion together with our
work in \S 2.2.3 (that determines the structure of $\euP_2^{\lceil
i/2\rceil}$) and determine the Galois module structure of $\euP_3^i$.
We need to express the image of $\mathcal{S}$ under the trace,
$\mbox{Tr}_{3,2}$, in terms of our $\euO_T[G]$-basis for
$\euP_2^{\lceil i/2\rceil}$.  This is the same as a determination of
the expression (in terms of Galois generators of $\euP_2^{\lceil
i/2q\rceil}$) for each valuation in (3.5).  First note under stable
ramification, $b_2 > 4e_0-2b_1$ so the structure of $\euP_2^{\lceil
i/2\rceil}$ is determined by the basis listed as Case $B$ in \S
2.2.3. However it is more convenient for us to use the basis listed as
Case $A$ in Appendix B. To translate between the two bases, note in
the elements $\overline{\alpha}$, $\overline{(\sigma +1)\alpha}$,
$\alpha$, $(\sigma +1)\alpha$ from \S 2.2.3 are referred to as
$\alpha$, $\rho$, $(\sigma^2 +1)\alpha$, $(\sigma +1)(\sigma^2
+1)\alpha$ in \S 3.1 and then in Appendix B. So row $B(1)$ in \S 2.2.3
corresponds with a pair of rows $A(7)$ and $A(8)$ in Appendix
B. Moreover $B(2)$ corresponds to rows $A(1)$ and $A(2)$, $B(3)$
corresponds to $A(3)$ and $A(4)$, and $B(4)$ corresponds to $A(5)$ and
$A(6)$.

There are four types of expression with valuation listed in (3.5). If
the valuation $a$ satisfies $a-(b_2-2b_1)\equiv 0\bmod 4$ then $a$ is
the valuation of a Galois translate $\rho$ where the difference between
$(\sigma \pm 1)\alpha$ and $\rho$ is an element $(\sigma +1)\mu\in K_0$
where $\mu$ is in the basis for $\euP_1^{i/4}$. Note  each such
$a$ corresponds with the appearance of $\mathcal{I}_2$ in the
$\euO_T[G]$ decomposition of $\euP_3^i$. Counting such $a$ one finds 
the same count as in $A(8)$. Note therefore  $A(7)$ counts the number of
$\mathcal{I}$ that are not mapped to under the trace, $\mbox{Tr}_{3,2}$, from
$\euP_3^i$.

Each valuation $a$ satisfying $a\equiv 0\bmod 4$ is the
valuation of $(\sigma^2 + 1)\alpha=(\sigma +1)\mu$ for some $\alpha$
in the basis for $\euP_2^{\lceil i/2\rceil}$ and $\mu$ in the basis
for $\euP_1^{i/4}\euP_0^{i/8}$. So each such $a$, corresponds with the
appearance of an $\mathcal{H}_2$. A count of such $a$ equals the count
in $A(2)$. Note $A(1)$ counts the number of $\mathcal{H}$ not
interacted with.  Each valuation $a$ satisfying $a-(b_2-2b_1)\equiv
2\bmod 4$ is the valuation of a Galois translate $\rho$ where the
difference between $(\sigma \pm 1)\alpha$ and $\rho$ is an element
$2\mu\in \euP_1^{i/4}$ where $(\sigma +1)\mu$ is in the basis for
$\euP_0^{i/8}$.  Each such $a$, therefore corresponds with the
appearance of an $\mathcal{M}_1$. The count of such $a$ is the same as
the count for $A(4)$. The number of $\mathcal{M}$ that appear in
$\euP_3^i$ is the same as the count for $A(3)$.  Finally each
valuation $a$ satisfying $a\equiv 2\bmod 4$ is the valuation of
$(\sigma^2 + 1)\alpha=2\mu$ for some $\alpha$ in the basis for
$\euP_2^{\lceil i/2\rceil}$.  Also $(\sigma +1)\mu$ is in the basis
for $\euP_0^{i/8}$, so each such $a$, therefore corresponds with the
appearance of an $\mathcal{L}_3$. The count of such $a$ is the same as
the count for $A(6)$. Again, $A(5)$ counts the number of $\mathcal{L}$
in $\euP_3^i$. 

Note the structure of $\euP_3^i$
under stable ramification is consistent with the structure of $\euP_3^i$
under unstable ramification so long as $b_2>4e_0-4b_1/3$.


\appendix
\section{The Modules}
In this section we introduce twenty--three indecomposable
$\mathbb{Z}_2[C_8]$-modules.  It is left to the interested reader to
translate our notation into Yakovlev's \cite{jakov:2}.  
\vspace*{3mm}

\noindent{\em Irreducibles}: Four of the $\mathbb{Z}_2[C_8]$-modules are
irreducible: $\mathcal{R}_0$, $\mathcal{R}_1$, $\mathcal{R}_2$, and
$\mathcal{R}_3$ where $\mathcal{R}_n:=\mathbb{Z}_2[\zeta_{2^n}]$,
$\zeta_{2^n}$ denotes a primitive $2^n$ root of unity, and $\sigma$
the generator of $C_8$ acts via multiplication by $\zeta_{2^n}$.

The other nineteen modules are `compounds'.  They are organized
according to fixed part -- those fixed by $\sigma^2$ are listed first,
followed by those fixed by $\sigma^4$, etc.
\vspace*{3mm}

\noindent{\em $\mathbb{Z}_2[C_2]$-modules}: Besides the two irreducibles
$\mathcal{R}_0$, $\mathcal{R}_1$, the group ring
$\mathbb{Z}_2[\sigma]/\langle\sigma^2\rangle$ is the only other
indecomposable module that is fixed by $\sigma^2$.  
\vspace*{3mm}

\noindent{\em Notation for `compounds'}: The group ring,
$\mathbb{Z}_2[\sigma]/\langle\sigma^2\rangle$, is made up of two
irreducibles. To make the relationships between irreducibles and their
`compounds' explicit, we will use diagrams like
$$\mathcal{R}_1\rightarrow 1\in \mathcal{R}_0$$ (instead of
$\mathbb{Z}_2[\sigma]/\langle\sigma^2\rangle$). These diagrams are to be
interpreted as follows: The number of $\mathbb{Z}_2[\sigma ]$-generators is
the number of irreducible modules that appear in the diagram. For
example, $\mathcal{R}_1\rightarrow 1\in \mathcal{R}_0$ means
two generators.  Let us call them $c$ and $d$. ({\em Think}: $c$
generates $\mathcal{R}_1$ while $d$ generates $\mathcal{R}_0$.)
Relations determine the module. If there is no `arrow' leaving an
irreducible $\mathcal{R}_i$, then the trace $\Phi_{2^i}(\sigma)$ maps
the generator to zero.  So $\Phi_{2^0}(\sigma)d=0$.  Note 
$\Phi_{2^i}(x)$ denotes the cyclotomic polynomial and
$x^8-1=\Phi_{2^0}(x)\cdot\Phi_{2^1}(x)\cdot\Phi_{2^2}(x)\cdot\Phi_{2^3}(x)$.
If there is an `arrow' leaving an irreducible $\mathcal{R}_i$
(pointing to an element), then the trace $\Phi_{2^i}(\sigma)$ maps the
generator to that element.  In this case $\Phi_{2^1}(\sigma)c =1\cdot
d$.
\vspace*{3mm}

\noindent{\em $\mathbb{Z}_2[C_4]$-modules}:
There are three indecomposable modules fixed by $\sigma^4$ (yet not
fixed by $\sigma^2$).  Notation for two other decomposable modules is
included as it will be needed to describe certain modules later (those not
fixed by $\sigma^4$). For three (of these five), the submodule
fixed by $\sigma^2$ is the group ring
$\mathbb{Z}_2[\sigma]/\langle\sigma^2\rangle$ (note how their diagams include
$\mathcal{R}_1\rightarrow 1\in \mathcal{R}_0$):
$$(\mathcal{G})\!: \mathcal{R}_2\rightarrow 1\in \mathcal{R}_1\rightarrow 1\in \mathcal{R}_0,\quad
(\mathcal{H})\!:\!\!\!\!
\begin{array}{lrr}
\mathcal{R}_2 & & \\
 &\searrow &  \\
\mathcal{R}_1 & \rightarrow&1\in \mathcal{R}_0
\end{array}\!\!\!,\quad
(\mathcal{I})\!: \mathcal{R}_2\oplus (\mathcal{R}_1\rightarrow 1\in
\mathcal{R}_0).$$ 
Denote the three generators by $b, c, d$.  ({\em
Think}: generating $\mathcal{R}_2, \mathcal{R}_1, \mathcal{R}_0$,
respectively.)  Recall $(\sigma-1)d=0$ while $(\sigma +1)c=d$.  In
$\mathcal{G}$, we have $\Phi_{2^2}(\sigma)b=1\cdot c$. So
$\mathcal{G}$ is the group ring
$\mathbb{Z}_2[\sigma]/\langle\sigma^4\rangle$. In $\mathcal{H}$, we have
$\Phi_{2^2}(\sigma)b=1\cdot d$.  While in $\mathcal{I}$,
$\Phi_{2^2}(\sigma)b=0$.  

For two (of these five), the submodule fixed
by $\sigma^2$ is the maximal order of
$\mathbb{Z}_2[\sigma]/\langle\sigma^2\rangle$ (note how
$\mathcal{R}_1\oplus\mathcal{R}_0$ appears):
$$(\mathcal{L}): 
\begin{array}{rrl}
\mathcal{R}_2& \rightarrow & 
\left \{ \begin{array}{l} 1\in \mathcal{R}_1 \\ \oplus  \\ 1\in \mathcal{R}_0 \end{array} \right .
\end{array},\quad
(\mathcal{M}): \mathcal{R}_2\oplus \mathcal{R}_1\oplus
\mathcal{R}_0.$$ 
Denote the three generators by $b, c, d$ where
$(\sigma-1)d=0$ and $(\sigma +1)c=0$.  In $\mathcal{L}$, we have
$\Phi_{2^2}(\sigma)b=1\cdot c+1\cdot d$.  In $\mathcal{M}$, we have
$\Phi_{2^2}(\sigma)b=0$.  So $\mathcal{M}$ is the maximal order of
$\mathbb{Z}_2[\sigma]/\langle\sigma^4\rangle$.
\vspace*{3mm}

\noindent{\em $\mathbb{Z}_2[C_8]$-modules}: The remaining fifteen
indecomposable modules can now be listed. They are collected according
to submodule fixed by $\sigma^4$.

\paragraph{\em Fixed part $\mathcal{G}$}
$$
\begin{array}{ll}
(\mathcal{G}_1): 
\mathcal{R}_3\rightarrow 1\in \mathcal{R}_2\rightarrow 1\in \mathcal{R}_1\rightarrow 1\in \mathcal{R}_0 &
(\mathcal{G}_3):
\begin{array}{rrrrr}\mathcal{R}_3& & & & \\
& \searrow& & & \\ \mathcal{R}_2& \rightarrow & 1\in \mathcal{R}_1&
\rightarrow & 1\in \mathcal{R}_0
\end{array} 
\vspace{4mm} \\
(\mathcal{G}_2): 
\mathcal{R}_3\rightarrow \lambda\in \mathcal{R}_2\rightarrow 1\in \mathcal{R}_1\rightarrow 1\in \mathcal{R}_0 &
(\mathcal{G}_4):
\begin{array}{rrrrr}& &\mathcal{R}_3 & & \\
& & &\searrow & \\
\mathcal{R}_2& \rightarrow & 1\in \mathcal{R}_1& \rightarrow & 1\in \mathcal{R}_0
\end{array} 
\end{array}$$
Call the generators $a, b, c, d$, where the $\mathbb{Z}_2[\sigma ]$-relations
among $b, c, d$ are as in $\mathcal{G}$.  In $\mathcal{G}_1$, we have
$\Phi_{2^3}(\sigma)a=1\cdot b$.  So $\mathcal{G}_1$ is the group ring
$\mathbb{Z}_2[\sigma]$.  In $\mathcal{G}_2$, we have
$\Phi_{2^3}(\sigma)a=\lambda\cdot b$ where $\lambda=\sigma-1$.  In
$\mathcal{G}_3$, $\Phi_{2^3}(\sigma)a=1\cdot c$.  In $\mathcal{G}_4$,
$\Phi_{2^3}(\sigma)a=1\cdot d$.

\paragraph{\em Fixed part $\mathcal{H}$}
$$\begin{array}{cc}
(\mathcal{H}_1): 
\begin{array}{rrl}
\mathcal{R}_3& \rightarrow & 
\left\{
\begin{array}{lrr}
\lambda\in \mathcal{R}_2 & & \\
\oplus &\searrow &  \\
1\in \mathcal{R}_1 & \rightarrow&1\in \mathcal{R}_0
\end{array} \right.
\end{array}
&
(\mathcal{H}_2):
\begin{array}{rrr}
\mathcal{R}_3 & & \\
 &\searrow & \\
\mathcal{R}_2& \rightarrow & 1\in \mathcal{R}_0 \\
 &\nearrow & \\
\mathcal{R}_1 & &
\end{array} 
\end{array}$$
Call the generators $a, b, c, d$, where the $\mathbb{Z}_2[\sigma
]$-relationships among $b, c, d$ are as in $\mathcal{H}$.  In
$\mathcal{H}_1$, $\Phi_{2^3}(\sigma)a=\lambda \cdot 1\cdot b+1\cdot
c$.  In $\mathcal{H}_2$, $\Phi_{2^3}(\sigma)a=d$.
\paragraph{\em Fixed part $\mathcal{I}$}
$$\begin{array}{ll}
(\mathcal{I}_1):
\begin{array}{rrl}
\mathcal{R}_3& \rightarrow & 
\left\{
\begin{array}{lrr}
1\in \mathcal{R}_2 & & \\
\oplus & &  \\
1\in \mathcal{R}_1 & \rightarrow&1\in \mathcal{R}_0
\end{array} \right.
\end{array} & 
(\mathcal{I}_2):
\begin{array}{rl}
\begin{array}{rr}
& \\
\mathcal{R}_3 &\rightarrow \\
\mathcal{R}_1& \rightarrow\end{array}\!\!
& \!\!\left \{
\begin{array}{l}
1\in \mathcal{R}_2 \\
\oplus  \\
1\in \mathcal{R}_0
\end{array}\right . 
\end{array} 
\end{array}$$
Each module is generated by $a, b, c, d$, where the $\mathbb{Z}_2[\sigma
]$-relationships among $b, c, d$ are as in $\mathcal{I}$. 
In $\mathcal{I}_1$, $\Phi_{2^3}(\sigma)a=
1\cdot b+1\cdot c$.  In $\mathcal{I}_2$, $\Phi_{2^3}(\sigma)a=1\cdot b+1\cdot d$.
\paragraph{\em Fixed part $\mathcal{L}$ or $\mathcal{M}$}
$$\begin{array}{ll} (\mathcal{L}_1):
\begin{array}{rrrrl}
\mathcal{R}_3& \rightarrow & 1\in \mathcal{R}_2& \rightarrow & 
\left \{ \begin{array}{l} 1\in \mathcal{R}_1 \\ \oplus  \\ 1\in \mathcal{R}_0 \end{array} \right .
\end{array} &
(\mathcal{L}_3):
\begin{array}{rrl}
 \mathcal{R}_3 & & \\ & \searrow& \\ \mathcal{R}_2 &\rightarrow& \left
\{ \begin{array}{l} 1\in \mathcal{R}_1 \\ \oplus \\ 1\in
\mathcal{R}_0 \end{array} \right .  \end{array}
\vspace*{4mm}\\ 
(\mathcal{L}_2):
\begin{array}{rrrrl}
\mathcal{R}_3& \rightarrow & \lambda\in \mathcal{R}_2& \rightarrow &
\left \{ \begin{array}{l} 1\in \mathcal{R}_1 \\ \oplus \\ 1\in
\mathcal{R}_0 \end{array} \right . \end{array} 
& (\mathcal{M}_1):
\begin{array}{rll}
\mathcal{R}_3& \rightarrow & 
\left \{ \begin{array}{l}
\lambda\in \mathcal{R}_2   \\
\oplus   \\
1\in \mathcal{R}_1   \\
\oplus   \\
1\in \mathcal{R}_0\end{array} \right .
\end{array}
\end{array}$$
The generators are
$a, b, c, d$, where the
$\mathbb{Z}_2[\sigma
]$-relationships among $b, c, d$ are as in $\mathcal{L}$ or
$\mathcal{M}$ respectively.  In $\mathcal{L}_1$, $\Phi_{2^3}(\sigma)a= b$.
In $\mathcal{L}_2$, $\Phi_{2^3}(\sigma)a=\lambda \cdot b$.  In
$\mathcal{L}_3$, $\Phi_{2^3}(\sigma)a=1\cdot c+1\cdot d$.  In $\mathcal{M}_1$,
$\Phi_{2^3}(\sigma)a=1\cdot b+1\cdot c+1\cdot d$.

\paragraph{\em Hybrids of $\mathcal{H}_1$}
The next three modules result from the
linking of an $\mathcal{H}_1$ with either another $\mathcal{R}_3$, or with a
$\mathcal{G}$, or with a $\mathcal{L}$.

$$
(\mathcal{H}_{1,2}):
\begin{array}{rrl}
& &\quad
\begin{array}{lrl}
\quad\;\; \mathcal{R}_3 & & \\ 
&\searrow & \end{array} \\ 
\mathcal{R}_3& \rightarrow & 
\left \{
\begin{array}{lrl}
1\in \mathcal{R}_1 & \rightarrow&1\in \mathcal{R}_0 \\
\oplus &\nearrow &  \\
\lambda\in \mathcal{R}_2 & & 
\end{array} \right .
\end{array} $$
This module is generated by $a_1, a_2, b, c, d$ with the
$\mathbb{Z}_2[\sigma]$-relationships among $b, c, d$ as in $\mathcal{H}$,
while $\Phi_{2^3}(\sigma)a_1=\lambda\cdot b+1\cdot c$ and $\Phi_{2^3}(\sigma)a_2=d$.  If
$\Phi_{2^3}(\sigma)a_1=0$, $\mathcal{H}_2$ would decompose off. If
$\Phi_{2^3}(\sigma)a_2=0$, $\mathcal{H}_1$  would decompose off. It is
a mixture of $\mathcal{H}_1$ and $\mathcal{H}_2$, hence the
name.

$$(\mathcal{H}_1\mathcal{G}):
\begin{array}{rrl}
\mathcal{R}_3& \rightarrow & 
	\left\{
		\begin{array}{lrl}
		1\in \mathcal{R}_1 & \rightarrow &1\in \mathcal{R}_0 \\
		\oplus & &  \\
		\lambda\in \mathcal{R}_2 &\rightarrow & \oplus 
		\end{array}
	\right . \\ 
& & \\
\mathcal{R}_2&\rightarrow & \quad
\begin{array}{lrl}
1\in \mathcal{R}_1 &\rightarrow & 1\in \mathcal{R}_0 
\end{array}
\end{array}
$$
This module is generated by $a_1, b_1, c_1, d_1$ and $b_2, c_2,
d_2$. The $\mathbb{Z}_2[\sigma]$-relationships among $b_2, c_2, d_2$ are as
in $\mathcal{G}$.  The $\mathbb{Z}_2[\sigma]$-relationships among $a_1, c_1,
d_1$ are as in $\mathcal{H}_1$ with $(\sigma^2+1)b_1=1\cdot d_1+1\cdot d_2$.

$$(\mathcal{H}_1\mathcal{L}):
\begin{array}{rrl}
\mathcal{R}_3& \rightarrow & 
\left\{
\begin{array}{lrl}
1\in \mathcal{R}_1 & \rightarrow& 1\in \mathcal{R}_0 \\
\oplus & &  \\
\lambda\in \mathcal{R}_2 &\rightarrow &\oplus 
\end{array} \right .
\\  
& & \\
& & 
\quad\;\; \begin{array}{lrl}
 \mathcal{R}_2 &\rightarrow & 
\left \{\begin{array}{l} 1\in \mathcal{R}_0 \\ 
\oplus \\
1\in \mathcal{R}_1\end{array} \right .\end{array}
\end{array}
$$
This module is generated by $a_1, b_1, c_1, d_1$ and $b_2, c_2,
d_2$. The $\mathbb{Z}_2[\sigma]$-relationships among $b_2, c_2, d_2$ are
as in $\mathcal{L}$.  The $\mathbb{Z}_2[\sigma]$-relationships among $a_1,
c_1, d_1$ are as in $\mathcal{H}_1$ with
$(\sigma^2+1)b_1=1\cdot d_1+(1\cdot c_2+1\cdot d_2)$.

\section{The Bases by Case, $A$ through $H$}
From \S 3.4, we inherit sequences of elements ordered in terms of
increasing valuation (for Case $A$, we have $\ldots \rho ,
\overline{2\rho}, (\sigma^2+1)\alpha, \overline{2(\sigma^2+1)\alpha},
2\alpha , \overline{4\alpha}, (\sigma +1)(\sigma^2+1)\alpha,
\overline{2(\sigma +1)(\sigma^2+1)\alpha}, 2\rho, \ldots $).
Following \S 2.2.3, we are interested in those elements `in view'
({\em i.e.} with valuation in $i, i+1, \ldots , i+v_3(2)-1$). As we
vary $m$ the `view' changes.  Indeed, for each case, there are eight
views (eight sets). They are listed below.  Recall from \S 2.2.3 
it is easy to determine the subscripts $m$ associated with a
particular `view'.  For example, the elements in $A(2)$ appear for $i\leq v_3\left
(\overline{(\sigma+1)(\sigma^2+1)\alpha}\right )$ and $v_3\left
((\sigma+1)(\sigma^2+1)\alpha\right)\leq 8e_0+i-1$.  In other words,
$\lceil(i+b_3-4b_1-4b_2)/8\rceil \leq m \leq \lceil
(i+8e_0-4b_1-4b_2)/8 \rceil -1$.
  

{\small
\begin{gather}
\mbox{\bf Case A} \notag\\ \rho,\; \overline{2\rho},\; (\sigma
^2+1)\alpha,\; \overline{2(\sigma ^2+1)\alpha},\; 2\alpha,\;
\overline{4\alpha},\; (\sigma+1)(\sigma^2+1)\alpha,\;
\overline{2(\sigma+1)(\sigma^2+1)\alpha} \tag{$1$}\\
\overline{(\sigma+1)(\sigma^2+1)\alpha},\; \rho,\; \overline{2\rho},\;
(\sigma ^2+1)\alpha,\; \overline{2(\sigma ^2+1)\alpha},\; 2\alpha,\;
\overline{4\alpha},\; (\sigma+1)(\sigma^2+1)\alpha \tag{$2$}\\
\frac{1}{2}(\sigma+1)(\sigma^2+1)\alpha,\;\overline{(\sigma+1)(\sigma^2+1)\alpha},\;
\rho,\; \overline{2\rho},\; (\sigma ^2+1)\alpha,\; \overline{2(\sigma
^2+1)\alpha},\; 2\alpha,\; \overline{4\alpha} \tag{$3$}\\
\overline{2\alpha},\;\frac{1}{2}(\sigma+1)(\sigma^2+1)\alpha,\;\overline{(\sigma+1)(\sigma^2+1)\alpha},\;
\rho,\; \overline{2\rho},\; (\sigma ^2+1)\alpha,\; \overline{2(\sigma
^2+1)\alpha},\; 2\alpha \tag{$4$}\\
\alpha,\;\overline{2\alpha},\;\frac{1}{2}(\sigma+1)(\sigma^2+1)\alpha,\;\overline{(\sigma+1)(\sigma^2+1)\alpha},\;
\rho,\; \overline{2\rho},\; (\sigma ^2+1)\alpha,\; \overline{2(\sigma
^2+1)\alpha} \tag{$5$}\\ \overline{(\sigma ^2+1)\alpha},\;
\alpha,\;\overline{2\alpha},\;\frac{1}{2}(\sigma+1)(\sigma^2+1)\alpha,\;\overline{(\sigma+1)(\sigma^2+1)\alpha},\;
\rho,\; \overline{2\rho},\; (\sigma ^2+1)\alpha \tag{$6$}\\
\frac{1}{2}(\sigma ^2+1)\alpha,\; \overline{(\sigma ^2+1)\alpha},\;
\alpha,\;\overline{2\alpha},\;\frac{1}{2}(\sigma+1)(\sigma^2+1)\alpha,\;\overline{(\sigma+1)(\sigma^2+1)\alpha},\;
\rho,\; \overline{2\rho} \tag{$7$}\\ \overline{\rho},\;
\frac{1}{2}(\sigma ^2+1)\alpha,\; \overline{(\sigma ^2+1)\alpha},\;
\alpha,\;\overline{2\alpha},\;\frac{1}{2}(\sigma+1)(\sigma^2+1)\alpha,\;\overline{(\sigma+1)(\sigma^2+1)\alpha},\;\rho
\tag{$8$}
\end{gather}
\begin{gather}
\mbox{\bf Case B} \notag\\
\rho,\; \overline{2\rho},\; (\sigma ^2+1)\alpha,\; \overline{2(\sigma ^2+1)\alpha},\; 2\alpha,\;  (\sigma+1)(\sigma^2+1)\alpha,\; \overline{4\alpha},\; \overline{2(\sigma+1)(\sigma^2+1)\alpha}  \tag{$1$}\\
\overline{(\sigma+1)(\sigma^2+1)\alpha},\;\rho,\; \overline{2\rho},\; (\sigma ^2+1)\alpha,\; \overline{2(\sigma ^2+1)\alpha},\; 2\alpha,\;  (\sigma+1)(\sigma^2+1)\alpha,\; \overline{4\alpha}   \tag{$2$}\\
\overline{2\alpha},\; \overline{(\sigma+1)(\sigma^2+1)\alpha},\;\rho,\; \overline{2\rho},\; (\sigma ^2+1)\alpha,\; \overline{2(\sigma ^2+1)\alpha},\; 2\alpha,\;  (\sigma+1)(\sigma^2+1)\alpha   \tag{$3$}\\
\frac{1}{2}(\sigma+1)(\sigma^2+1)\alpha,\; \overline{2\alpha},\; \overline{(\sigma+1)(\sigma^2+1)\alpha},\;\rho,\; \overline{2\rho},\; (\sigma ^2+1)\alpha,\; \overline{2(\sigma ^2+1)\alpha},\; 2\alpha    \tag{$4$}\\
\alpha,\;  \frac{1}{2}(\sigma+1)(\sigma^2+1)\alpha,\; \overline{2\alpha},\; \overline{(\sigma+1)(\sigma^2+1)\alpha},\;\rho,\; \overline{2\rho},\; (\sigma ^2+1)\alpha,\; \overline{2(\sigma ^2+1)\alpha}   \tag{$5$}\\
\overline{(\sigma ^2+1)\alpha},\;\alpha,\;  \frac{1}{2}(\sigma+1)(\sigma^2+1)\alpha,\; \overline{2\alpha},\; \overline{(\sigma+1)(\sigma^2+1)\alpha},\;\rho,\; \overline{2\rho},\; (\sigma ^2+1)\alpha   \tag{$6$}\\
\frac{1}{2}(\sigma ^2+1)\alpha,\; \overline{(\sigma ^2+1)\alpha},\;\alpha,\;  \frac{1}{2}(\sigma+1)(\sigma^2+1)\alpha,\; \overline{2\alpha},\; \overline{(\sigma+1)(\sigma^2+1)\alpha},\;\rho,\; \overline{2\rho}   \tag{$7$}\\
 \overline{\rho},\;\frac{1}{2}(\sigma ^2+1)\alpha,\; \overline{(\sigma ^2+1)\alpha},\;\alpha,\;  \frac{1}{2}(\sigma+1)(\sigma^2+1)\alpha,\; \overline{2\alpha},\; \overline{(\sigma+1)(\sigma^2+1)\alpha},\;\rho    \tag{$8$}
\end{gather}
\begin{gather}
\mbox{\bf Case C} \notag\\
\rho,\; \overline{2\rho},\; (\sigma ^2+1)\alpha,\; \overline{2(\sigma ^2+1)\alpha},\;(\sigma+1)(\sigma^2+1)\alpha,\; 2\alpha,\; \overline{2(\sigma+1)(\sigma^2+1)\alpha}, \;   \overline{4\alpha} \tag{$1$}\\
\overline{2\alpha},\;\rho,\; \overline{2\rho},\; (\sigma ^2+1)\alpha,\; \overline{2(\sigma ^2+1)\alpha},\;(\sigma+1)(\sigma^2+1)\alpha,\; 2\alpha,\; \overline{2(\sigma+1)(\sigma^2+1)\alpha}    \tag{$2$}\\
\overline{(\sigma+1)(\sigma^2+1)\alpha}, \; \overline{2\alpha},\;\rho,\; \overline{2\rho},\; (\sigma ^2+1)\alpha,\; \overline{2(\sigma ^2+1)\alpha},\;(\sigma+1)(\sigma^2+1)\alpha,\; 2\alpha    \tag{$3$}\\
\alpha,\; \overline{(\sigma+1)(\sigma^2+1)\alpha}, \; \overline{2\alpha},\;\rho,\; \overline{2\rho},\; (\sigma ^2+1)\alpha,\; \overline{2(\sigma ^2+1)\alpha},\;(\sigma+1)(\sigma^2+1)\alpha    \tag{$4$}\\
\frac{1}{2}(\sigma+1)(\sigma^2+1)\alpha,\;\alpha,\; \overline{(\sigma+1)(\sigma^2+1)\alpha}, \; \overline{2\alpha},\;\rho,\; \overline{2\rho},\; (\sigma ^2+1)\alpha,\; \overline{2(\sigma ^2+1)\alpha}    \tag{$5$}\\
\overline{(\sigma ^2+1)\alpha},\;\frac{1}{2}(\sigma+1)(\sigma^2+1)\alpha,\;\alpha,\; \overline{(\sigma+1)(\sigma^2+1)\alpha}, \; \overline{2\alpha},\;\rho,\; \overline{2\rho},\; (\sigma ^2+1)\alpha     \tag{$6$}\\
\frac{1}{2}(\sigma ^2+1)\alpha,\; \overline{(\sigma ^2+1)\alpha},\;\frac{1}{2}(\sigma+1)(\sigma^2+1)\alpha,\;\alpha,\; \overline{(\sigma+1)(\sigma^2+1)\alpha}, \; \overline{2\alpha},\;\rho,\; \overline{2\rho}    \tag{$7$}\\
\overline{\rho},\; \frac{1}{2}(\sigma ^2+1)\alpha,\; \overline{(\sigma ^2+1)\alpha},\;\frac{1}{2}(\sigma+1)(\sigma^2+1)\alpha,\;\alpha,\; \overline{(\sigma+1)(\sigma^2+1)\alpha}, \; \overline{2\alpha},\;\rho     \tag{$8$}
\end{gather}

\begin{gather}
\mbox{\bf Case D} \notag\\
\rho,\; (\sigma ^2+1)\alpha,\;\overline{2\rho},\;  \overline{2(\sigma ^2+1)\alpha},\;(\sigma+1)(\sigma^2+1)\alpha,\; 2\alpha,\; \overline{2(\sigma+1)(\sigma^2+1)\alpha}, \;   \overline{4\alpha} \tag{$1$}\\
\overline{2\alpha}, \;\rho,\; (\sigma ^2+1)\alpha,\;\overline{2\rho},\;  \overline{2(\sigma ^2+1)\alpha},\;(\sigma+1)(\sigma^2+1)\alpha,\; 2\alpha,\; \overline{2(\sigma+1)(\sigma^2+1)\alpha}    \tag{$2$}\\
\overline{(\sigma+1)(\sigma^2+1)\alpha}, \;\overline{2\alpha}, \;\rho,\; (\sigma ^2+1)\alpha,\;\overline{2\rho},\;  \overline{2(\sigma ^2+1)\alpha},\;(\sigma+1)(\sigma^2+1)\alpha,\; 2\alpha     \tag{$3$}\\
\alpha,\;\overline{(\sigma+1)(\sigma^2+1)\alpha}, \;\overline{2\alpha}, \;\rho,\; (\sigma ^2+1)\alpha,\;\overline{2\rho},\;  \overline{2(\sigma ^2+1)\alpha},\;(\sigma+1)(\sigma^2+1)\alpha      \tag{$4$}\\
\frac{1}{2}(\sigma+1)(\sigma^2+1)\alpha,\;\alpha,\;\overline{(\sigma+1)(\sigma^2+1)\alpha}, \;\overline{2\alpha}, \;\rho,\; (\sigma ^2+1)\alpha,\;\overline{2\rho},\;  \overline{2(\sigma ^2+1)\alpha}      \tag{$5$}\\
 \overline{(\sigma ^2+1)\alpha},\; \frac{1}{2}(\sigma+1)(\sigma^2+1)\alpha,\;\alpha,\;\overline{(\sigma+1)(\sigma^2+1)\alpha}, \;\overline{2\alpha}, \;\rho,\; (\sigma ^2+1)\alpha,\;\overline{2\rho}      \tag{$6$}\\
\overline{\rho},\; \overline{(\sigma ^2+1)\alpha},\; \frac{1}{2}(\sigma+1)(\sigma^2+1)\alpha,\;\alpha,\;\overline{(\sigma+1)(\sigma^2+1)\alpha}, \;\overline{2\alpha}, \;\rho,\; (\sigma ^2+1)\alpha      \tag{$7$}\\
\frac{1}{2}(\sigma ^2+1)\alpha,\;\overline{\rho},\; \overline{(\sigma ^2+1)\alpha},\; \frac{1}{2}(\sigma+1)(\sigma^2+1)\alpha,\;\alpha,\;\overline{(\sigma+1)(\sigma^2+1)\alpha}, \;\overline{2\alpha}, \;\rho       \tag{$8$}
\end{gather}

\begin{gather}
\mbox{\bf Case E} \notag\\
\overline{2\alpha},\; \overline{2\rho},\;(\sigma ^2+1)\alpha,\;(\sigma+1)(\sigma^2+1)\alpha,\; \overline{2(\sigma ^2+1)\alpha},\;\overline{2(\sigma+1)(\sigma^2+1)\alpha}, \;2\alpha,\;2\rho      \tag{$1$}\\
\rho,\; \overline{2\alpha},\;\overline{2\rho},\; (\sigma ^2+1)\alpha,\;  (\sigma+1)(\sigma^2+1)\alpha,\; \overline{2(\sigma ^2+1)\alpha},\; \overline{2(\sigma+1)(\sigma^2+1)\alpha}, \; 2\alpha   \tag{$2$}\\
\alpha,\;\rho,\; \overline{2\alpha},\;\overline{2\rho},\; (\sigma ^2+1)\alpha,\;  (\sigma+1)(\sigma^2+1)\alpha,\; \overline{2(\sigma ^2+1)\alpha},\; \overline{2(\sigma+1)(\sigma^2+1)\alpha}    \tag{$3$}\\
\overline{(\sigma+1)(\sigma^2+1)\alpha}, \;\alpha,\;\rho,\; \overline{2\alpha},\;\overline{2\rho},\; (\sigma ^2+1)\alpha,\;  (\sigma+1)(\sigma^2+1)\alpha,\; \overline{2(\sigma ^2+1)\alpha}     \tag{$4$}\\
\overline{(\sigma ^2+1)\alpha},\;\overline{(\sigma+1)(\sigma^2+1)\alpha}, \;\alpha,\;\rho,\; \overline{2\alpha},\;\overline{2\rho},\; (\sigma ^2+1)\alpha,\;  (\sigma+1)(\sigma^2+1)\alpha     \tag{$5$}\\
\frac{1}{2}(\sigma+1)(\sigma^2+1)\alpha,\; \overline{(\sigma ^2+1)\alpha},\;\overline{(\sigma+1)(\sigma^2+1)\alpha}, \;\alpha,\;\rho,\; \overline{2\alpha},\;\overline{2\rho},\; (\sigma ^2+1)\alpha      \tag{$6$}\\
\frac{1}{2}(\sigma ^2+1)\alpha,\;\frac{1}{2}(\sigma+1)(\sigma^2+1)\alpha,\; \overline{(\sigma ^2+1)\alpha},\;\overline{(\sigma+1)(\sigma^2+1)\alpha}, \;\alpha,\;\rho,\; \overline{2\alpha},\;\overline{2\rho}      \tag{$7$}\\
\overline{\rho},\;\frac{1}{2}(\sigma ^2+1)\alpha,\;\frac{1}{2}(\sigma+1)(\sigma^2+1)\alpha,\; \overline{(\sigma ^2+1)\alpha},\;\overline{(\sigma+1)(\sigma^2+1)\alpha}, \;\alpha,\;\rho,\; \overline{2\alpha}      \tag{$8$}
\end{gather}

\begin{gather}
\mbox{\bf Case F} \notag\\
\overline{2\alpha},\;(\sigma ^2+1)\alpha,\; \overline{2\rho},\;(\sigma+1)(\sigma^2+1)\alpha,\;\overline{2(\sigma ^2+1)\alpha},\;\overline{2(\sigma+1)(\sigma^2+1)\alpha}, \;2\alpha,\;2\rho        \tag{$1$} \\
\rho,\; \overline{2\alpha},\;(\sigma ^2+1)\alpha,\;\overline{2\rho},\;   (\sigma+1)(\sigma^2+1)\alpha,\; \overline{2(\sigma ^2+1)\alpha},\; \overline{2(\sigma+1)(\sigma^2+1)\alpha}, \; 2\alpha   \tag{$2$}\\
\alpha,\;\rho,\; \overline{2\alpha},\;(\sigma ^2+1)\alpha,\;\overline{2\rho},\;   (\sigma+1)(\sigma^2+1)\alpha,\; \overline{2(\sigma ^2+1)\alpha},\; \overline{2(\sigma+1)(\sigma^2+1)\alpha}    \tag{$3$}\\
\overline{(\sigma+1)(\sigma^2+1)\alpha}, \;\alpha,\;\rho,\; \overline{2\alpha},\;(\sigma ^2+1)\alpha,\;\overline{2\rho},\;   (\sigma+1)(\sigma^2+1)\alpha,\; \overline{2(\sigma ^2+1)\alpha}    \tag{$4$}\\
\overline{(\sigma ^2+1)\alpha},\;\overline{(\sigma+1)(\sigma^2+1)\alpha}, \;\alpha,\;\rho,\; \overline{2\alpha},\;(\sigma ^2+1)\alpha,\;\overline{2\rho},\;   (\sigma+1)(\sigma^2+1)\alpha     \tag{$5$}\\
\frac{1}{2}(\sigma+1)(\sigma^2+1)\alpha,\;\overline{(\sigma ^2+1)\alpha},\;\overline{(\sigma+1)(\sigma^2+1)\alpha}, \;\alpha,\;\rho,\; \overline{2\alpha},\;(\sigma ^2+1)\alpha,\;\overline{2\rho}        \tag{$6$}\\
\overline{\rho},\;\frac{1}{2}(\sigma+1)(\sigma^2+1)\alpha,\;\overline{(\sigma ^2+1)\alpha},\;\overline{(\sigma+1)(\sigma^2+1)\alpha}, \;\alpha,\;\rho,\; \overline{2\alpha},\;(\sigma ^2+1)\alpha        \tag{$7$}\\
\frac{1}{2}(\sigma ^2+1)\alpha,\; \overline{\rho},\;\frac{1}{2}(\sigma+1)(\sigma^2+1)\alpha,\;\overline{(\sigma ^2+1)\alpha},\;\overline{(\sigma+1)(\sigma^2+1)\alpha}, \;\alpha,\;\rho,\; \overline{2\alpha}       \tag{$8$}
\end{gather}

\begin{gather}
\mbox{\bf Case G} \notag\\
(\sigma ^2+1)\alpha,\;\overline{2\alpha},\;(\sigma+1)(\sigma^2+1)\alpha,\; \overline{2\rho},\;\overline{2(\sigma ^2+1)\alpha},\;\overline{2(\sigma+1)(\sigma^2+1)\alpha}, \;2\alpha,\;2\rho    \tag{$1$}\\
\rho,\;(\sigma ^2+1)\alpha,\; \overline{2\alpha},\;(\sigma+1)(\sigma^2+1)\alpha,\; \overline{2\rho},\;\overline{2(\sigma ^2+1)\alpha},\; \overline{2(\sigma+1)(\sigma^2+1)\alpha}, \; 2\alpha   \tag{$2$}\\
\alpha,\;\rho,\;(\sigma ^2+1)\alpha,\; \overline{2\alpha},\;(\sigma+1)(\sigma^2+1)\alpha,\; \overline{2\rho},\;\overline{2(\sigma ^2+1)\alpha},\; \overline{2(\sigma+1)(\sigma^2+1)\alpha}   \tag{$3$}\\
\overline{(\sigma+1)(\sigma^2+1)\alpha}, \;\alpha,\;\rho,\;(\sigma ^2+1)\alpha,\; \overline{2\alpha},\;(\sigma+1)(\sigma^2+1)\alpha,\; \overline{2\rho},\;\overline{2(\sigma ^2+1)\alpha}    \tag{$4$}\\
\overline{(\sigma ^2+1)\alpha},\;\overline{(\sigma+1)(\sigma^2+1)\alpha}, \;\alpha,\;\rho,\;(\sigma ^2+1)\alpha,\; \overline{2\alpha},\;(\sigma+1)(\sigma^2+1)\alpha,\; \overline{2\rho}   \tag{$5$}\\
\overline{\rho},\;\overline{(\sigma ^2+1)\alpha},\;\overline{(\sigma+1)(\sigma^2+1)\alpha}, \;\alpha,\;\rho,\;(\sigma ^2+1)\alpha,\; \overline{2\alpha},\;(\sigma+1)(\sigma^2+1)\alpha    \tag{$6$}\\
\frac{1}{2}(\sigma+1)(\sigma^2+1)\alpha,\; \overline{\rho},\;\overline{(\sigma ^2+1)\alpha},\;\overline{(\sigma+1)(\sigma^2+1)\alpha}, \;\alpha,\;\rho,\;(\sigma ^2+1)\alpha,\; \overline{2\alpha}   \tag{$7$}\\
\overline{\alpha},\;\frac{1}{2}(\sigma+1)(\sigma^2+1)\alpha,\; \overline{\rho},\;\overline{(\sigma ^2+1)\alpha},\;\overline{(\sigma+1)(\sigma^2+1)\alpha}, \;\alpha,\;\rho,\;(\sigma ^2+1)\alpha    \tag{$8$}
\end{gather}

\begin{gather}
\mbox{\bf Case H} \notag\\
(\sigma ^2+1)\alpha,\;(\sigma+1)(\sigma^2+1)\alpha,\;\overline{2\alpha},\;\overline{2\rho},\; \overline{2(\sigma ^2+1)\alpha},\;\overline{2(\sigma+1)(\sigma^2+1)\alpha}, \;2\alpha,\;2\rho  \tag{$1$}\\
\rho,\;(\sigma ^2+1)\alpha,\; (\sigma+1)(\sigma^2+1)\alpha,\;\overline{2\alpha},\; \overline{2\rho},\;\overline{2(\sigma ^2+1)\alpha},\; \overline{2(\sigma+1)(\sigma^2+1)\alpha}, \; 2\alpha   \tag{$2$}\\
\alpha,\;\rho,\;(\sigma ^2+1)\alpha,\; (\sigma+1)(\sigma^2+1)\alpha,\;\overline{2\alpha},\; \overline{2\rho},\;\overline{2(\sigma ^2+1)\alpha},\; \overline{2(\sigma+1)(\sigma^2+1)\alpha}    \tag{$3$}\\
\overline{(\sigma+1)(\sigma^2+1)\alpha}, \;\alpha,\;\rho,\;(\sigma ^2+1)\alpha,\; (\sigma+1)(\sigma^2+1)\alpha,\;\overline{2\alpha},\; \overline{2\rho},\;\overline{2(\sigma ^2+1)\alpha}    \tag{$4$}\\
\overline{(\sigma ^2+1)\alpha},\;\overline{(\sigma+1)(\sigma^2+1)\alpha}, \;\alpha,\;\rho,\;(\sigma ^2+1)\alpha,\; (\sigma+1)(\sigma^2+1)\alpha,\;\overline{2\alpha},\; \overline{2\rho}    \tag{$5$}\\
\overline{\rho},\; \overline{(\sigma ^2+1)\alpha},\;\overline{(\sigma+1)(\sigma^2+1)\alpha}, \;\alpha,\;\rho,\;(\sigma ^2+1)\alpha,\; (\sigma+1)(\sigma^2+1)\alpha,\;\overline{2\alpha}    \tag{$6$}\\
\overline{\alpha},\;\overline{\rho},\; \overline{(\sigma ^2+1)\alpha},\;\overline{(\sigma+1)(\sigma^2+1)\alpha}, \;\alpha,\;\rho,\;(\sigma ^2+1)\alpha,\; (\sigma+1)(\sigma^2+1)\alpha   \tag{$7$}\\
\frac{1}{2}(\sigma+1)(\sigma^2+1)\alpha,\;\overline{\alpha},\;\overline{\rho},\; \overline{(\sigma ^2+1)\alpha},\;\overline{(\sigma+1)(\sigma^2+1)\alpha}, \;\alpha,\;\rho,\;(\sigma ^2+1)\alpha    \tag{$8$}
\end{gather}
}
\bibliographystyle{alpha}
\bibliography{bib}

\begin{thebibliography}{RCVSM90}

\bibitem[BE02]{elder:byott}
N.~P. Byott and G.~G. Elder.
\newblock Biquadratic extensions with one break.
\newblock {\em Canad. Math. Bull.}, 45(2):168--179, 2002.

\bibitem[CR90]{curt}
C.~W. Curtis and I.~Reiner.
\newblock {\em {Methods of Representation Theory}}.
\newblock Wiley-Interscience, New York, 1990.

\bibitem[Die85]{diet}
E.~Dieterich.
\newblock Representation types of group rings over complete discrete valuation
  rings. {I}{I}.
\newblock In {\em Orders and their applications (Oberwolfach, 1984)}, pages
  112--125. Springer, Berlin, 1985.

\bibitem[Eld95]{elder:annals}
G.~G. Elder.
\newblock {Galois module structure of integers in wildly ramified cyclic
  extensions of degree $p^2$}.
\newblock {\em Ann. Inst. Fourier (Grenoble)}, 45(3):625--647, 1995.
\newblock {{\em errata ibid.} {\bf 48} (1998), no. 2, 609--610}.

\bibitem[Eld98]{elder:biquad}
G.~G. Elder.
\newblock {Galois module structure of ideals in wildly ramified biquadratic
  extensions}.
\newblock {\em Can. J. Math.}, 50(5):1007--1047, 1998.

\bibitem[Eld02]{elder:bord}
G.~G. Elder.
\newblock {On Galois structure of the integers in cyclic extensions of local
  number fields}.
\newblock {\em J. Th\'{e}or. Nombres Bordeaux}, 14(1):113--149, 2002.

\bibitem[EM94]{elder:jnt}
G.~G. Elder and M.~L. Madan.
\newblock {Galois module structure of integers in wildly ramified cyclic
  extensions}.
\newblock {\em J. Number Theory}, 47(2):138--174, 1994.

\bibitem[Fon71]{fontaine}
J.-M. Fontaine.
\newblock {Groupes de ramification et repr\'{e}sentations d'Artin}.
\newblock {\em Ann. Scient. \'{E}c. Norm. Sup.}, 4:337--392, 1971.

\bibitem[HKO98]{kling:ks}
P.~Hindman, L.~Klingler, and C.~J. Odenthal.
\newblock On the {K}rull-{S}chmidt-{A}zumaya theorem for integral group rings.
\newblock {\em Comm. Algebra}, 26(11):3743--3758, 1998.

\bibitem[Jak75]{jakov:2}
A.~V. Jakovlev.
\newblock {Classification of $2$-adic Representations of an Eighth-Order Cyclic
  Group}.
\newblock {\em Journal of Soviet Math.}, 3(5):654--680, 1975.

\bibitem[Miy87]{miyata}
Y.~Miyata.
\newblock Vertices of ideals of a $p$-adic number field.
\newblock {\em Illinois J. Math.}, 31(2):185--199, 1987.

\bibitem[Miy95]{miyata:2}
Y.~Miyata.
\newblock On the {G}alois module structure of ideals and rings of all integers
  of {$p$}-adic number fields.
\newblock {\em J. Algebra}, 177(3):627--646, 1995.

\bibitem[Naz61]{naza}
L.~A. Nazarova.
\newblock Integral representations of klein's four-group.
\newblock {\em Soviet Math. Dokl.}, 2:1304--1307, 1961.
\newblock English Translation.

\bibitem[Noe32]{emmy}
E.~Noether.
\newblock {Normalbasis bei K\"{o}rpern ohne h\"{o}here Verzweigung}.
\newblock {\em J. Reine Angew. Math.}, 167:147--152, 1932.

\bibitem[RCVSM90]{martha}
M.~Rzedowski-Calder\'{o}n, G.~Villa-Salvador, and M.~L. Madan.
\newblock {Galois module structure of rings of integers}.
\newblock {\em Math. Z.}, 204:401--424, 1990.

\bibitem[Ser79]{serre:local}
J-P. Serre.
\newblock {\em {Local fields}}.
\newblock Springer-Verlag, Berlin/Heidelberg/New York, 1979.

\bibitem[Ull70]{ullom}
S.~Ullom.
\newblock Integral normal bases in {G}alois extensions of local fields.
\newblock {\em Nagoya Math. J.}, 39:141--148, 1970.

\bibitem[Wie84]{wiegand}
R.~Wiegand.
\newblock Cancellation over commutative rings of dimension one and two.
\newblock {\em J. Algebra}, 88(2):438--459, 1984.

\bibitem[Wym69]{wyman}
B.~Wyman.
\newblock {Wildly ramified gamma extensions}.
\newblock {\em Am. J. Math.}, 91:135--152, 1969.

\end{thebibliography}
\end{document}